\documentclass[11pt]{article}

\usepackage{hyperref}
\usepackage[T1]{fontenc}
\usepackage[latin1]{inputenc}
\usepackage{amsmath,amsthm,amssymb}
\usepackage{mathabx}
\usepackage{xcolor}
\usepackage{authblk}
\usepackage[margin=1in]{geometry}
\usepackage{eufrak}
\usepackage{mathrsfs}
\usepackage{enumitem} 

\providecommand{\U}[1]{\protect\rule{.1in}{.1in}}

\newtheorem{prop}{Proposition}[section]
\newtheorem{cor}[prop]{Corollary}

\newtheorem{theo}[prop]{Theorem}

\newenvironment{remark}[1][Remark:\\]{\begin{trivlist}\item[\hskip \labelsep {\bfseries #1}]}{\end{trivlist}}

\newcommand{\vertiii}[1]{{\left\vert\kern-0.25ex\left\vert\kern-0.25ex\left\vert #1
    \right\vert\kern-0.25ex\right\vert\kern-0.25ex\right\vert}}

\def\tr{\mbox{\rm Tr}}

\newcommand{\EE}{\mathbb{E}}

\newcommand{\PP}{\mathbb{P}}

\newcommand{\RR}{\mathbb{R}}

\newcommand{\Aa}{ {\cal A }}
\newcommand{\Ba}{ {\cal B }}
\newcommand{\Ca}{ {\cal C }}
\newcommand{\Da}{ {\cal D }}

\def\As{{\mathscr A}}

\newcommand{\Ea}{ {\cal E }}

\newcommand{\Va}{ {\cal V }}

\newcommand{\Fa}{ {\cal F }}
\newcommand{\Ga}{ {\cal G }}

\newcommand{\Ma}{ {\cal M }}

\newcommand{\Ha}{ {\cal H }}

\newcommand{\Pa}{ {\cal P }}
\newcommand{\Za}{ {\cal Z }}

\newcommand{\Wa}{ {\cal W }}

\newcommand{\point}{\mbox{\LARGE .}}

\newcommand{\cqfd}{\hfill\blbx \\}
\def\blbx{\hbox{\vrule height 5pt width 5pt depth 0pt}\medskip}

\def \PP{\mathbb{P}}
\def \RR{\mathbb{R}}
\def \EE{\mathbb{E}}
\def \EE{\mathbb{E}}

\def \WW{\mathbb{W}}

\numberwithin{equation}{section}

\makeatletter
\DeclareRobustCommand\frownotimes{\mathbin{\mathpalette\frown@otimes\relax}}
\newcommand{\frown@otimes}[2]{
  \vbox{
    \ialign{##\cr
      \hidewidth$\m@th#1{}_\frown$\kern-\scriptspace\hidewidth\cr
      \noalign{\nointerlineskip\kern-.1pt}
      $\m@th#1\otimes$\cr
    }
  }
}
\makeatother

\newcommand\quotient[2]{
        \mathchoice
            {
                \text{\raise0ex\hbox{$#1$}/\lower0ex\hbox{$#2$}}%
            }
            {
                #1\,/\,#2
            }
            {
                #1\,/\,#2
            }
            {
                #1\,/\,#2
            }
    }

\begin{document}

\title{A variational approach to nonlinear and interacting diffusions}
\author[$1$]{M. Arnaudon}
\author[$2$]{P. Del Moral\thanks{P. Del Moral was supported in part by funding from the Chaire Stress Test, BNP Paribas SFTS and CMAP, Polytechnique Palaiseau, France}}
\affil[$1$]{{\small Institut de Math\'ematiques de Bordeaux (IMB), Bordeaux University, France}}
\affil[$2$]{{\small INRIA, Bordeaux Research Center \& CMAP, Polytechnique Palaiseau, France}}
\date{}

\maketitle

\begin{abstract}
The article presents a novel variational calculus to analyze the stability and the propagation of chaos properties of nonlinear and interacting diffusions.
This differential methodology combines gradient flow estimates with backward stochastic interpolations, Lyapunov linearization techniques as well as spectral theory. 
This framework applies to
a large class of stochastic models including non homogeneous diffusions, as well as stochastic processes evolving on 
differentiable manifolds, such as
constraint-type embedded  
manifolds on Euclidian spaces and manifolds equipped with some Riemannian metric.  We derive uniform as well as almost sure exponential contraction inequalities at the level of the nonlinear diffusion flow, yielding what seems to be the first result of this type for this class of models. Uniform propagation of chaos properties w.r.t. the time parameter are also provided.
 Illustrations are provided in the context of a class of gradient flow diffusions arising in fluid mechanics and granular media literature. The extended versions of these nonlinear Langevin-type diffusions on Riemannian manifolds are also discussed.\\



\end{abstract}

\section{Introduction}
\subsection{Description of the models}
We denote by $\Vert A\Vert_{2}:=\lambda_{\tiny max}(AA^{\prime})^{1/2}$, resp. $\Vert A\Vert_{F}=\tr(AA^{\prime})^{1/2}$ and $\rho(A)=\lambda_{\tiny max}((A+A^{\prime})/2)$  the spectral norm, the Frobenius norm, and the logarithmic norm of some matrix $A$, where $A^{\prime}$ stands for the transpose of $A$, and $\lambda_{\tiny max}(\point)$ the maximal eigenvalue.
With a slight abuse of notation, we denote by $I$ the identity $(d\times d)$-matrix, for any $d\geq 1$.    

Let $b_t$ be some time varying differentiable vector field with Jacobian matrix $\nabla b_t$ on $\RR^d$, for some parameter $d\geq 1$.
Consider the deterministic flow $t\in [s,\infty[\mapsto X_{s,t}(x)$
starting at $X_{s,s}(x)=x$ associated with the evolution equation 
  \begin{equation}\label{taylor-first-0}
\partial_t X_{s,t}(x)=b_t(X_{s,t}(x))\Longrightarrow \partial_t \nabla X_{s,t}(x)=\nabla X_{s,t}(x)~\nabla b_t(X_{s,t}(x))\quad \mbox{\rm with}\quad
\nabla X_{s,s}(x)=I
    \end{equation}
The r.h.s. equation is often called the first order variational equation associated with the flow $X_{s,t}(x)$ along the trajectory $X_{s,t}(x)$.
This equation plays a central role in the sensitivity analysis of nonlinear dynamical systems w.r.t. their initial conditions. 
For instance, the spectral norm of $\nabla X_{s,t}(x)$ can be estimated in terms of the logarithmic norm using the inequalities
    \begin{equation}\label{coppel-ineq}
-\int_s^t~\rho\left(-\nabla b_u(X_{s,u}(x))\right) ~du\leq \log{\Vert \nabla X_{s,t}(x)\Vert_2}\leq\int_s^t~\rho\left(\nabla b_u(X_{s,u}(x))\right) ~du
    \end{equation}
A proof of this assertion can be found in~\cite{coppel}, see also~\cite{martin} for extensions to Lipschitz functions on Banach spaces. 
Whenever $\rho\left(\nabla b_u(x)\right)\leq -\lambda$ for some $\lambda>0$, the r.h.s. estimate in (\ref{coppel-ineq}) readily implies the exponential stability estimate
  \begin{equation}\label{taylor-first}
\begin{array}{l}
\displaystyle X_{s,t}(x)-X_{s,t}(y)=\int_0^1~\langle\nabla X_{s,t}(\epsilon x+(1-\epsilon)y),(x-y)\rangle~d\epsilon\\
\\
\displaystyle \Longrightarrow \Vert X_{s,t}(x)-X_{s,t}(y)\Vert\leq e^{-\lambda(t-s)}~\Vert x-y\Vert
\end{array}
    \end{equation}
    
The linearization technique discussed above is often referred as the Lyapunov first or indirect method to analyze the stability of nonlinear dynamical systems. For a more thorough discussion on this subject we refer to the pioneering work by Lyapunov~\cite{lyapunov}, as well as to chapter 4 in the more recent monograph by Khalil~\cite{khalil}. 

The main objective of this article is to extend these results to 
nonlinear diffusions and their mean field particle interpretations on Euclidian as well as on 
differentiable manifolds.  The differential analysis of conventional diffusions w.r.t. initial conditions is also one of the stepping stones of 
Bismut and Malliavin  calculus. This framework is mainly designed to study the existence and the properties of smooth probability densities in terms of the differential properties of the diffusion semigroup. For a more thorough discussion on this subject we refer to~\cite{carverhill,norris}, and references therein.

The relevant mathematical apparatus for the description and the variational analysis of stochastic processes on manifolds being technically more sophisticated than conventional differential calculus, this introduction only discusses nonlinear and interacting diffusions on Euclidian spaces.   The extended versions of these models on Riemannian manifolds are discussed in some details
 in section~\ref{sec-manifold}, as well as in section~\ref{subsec4.3}.

Let $\Pa_2(\RR^d)$ be the set of Borel probability measures  on $\RR^d$ with finite second absolute moment,
equipped  with the $2$-Wasserstein distance given by
$$
\WW_2(\eta,\mu)=\inf\EE(\Vert X-Y\Vert^2)^{1/2}$$
In the above display, the infimum is taken over all pairs of random variables $(X,Y)$ with respective distributions $\eta$ and $\mu\in \Pa_2(\RR^d)$; and $\Vert X-Y\Vert$ stands for the Euclidian distance between $X$ and $Y$ on the product space $\RR^d$.

Also let  $b_t$ and $\sigma_{t}$ be differentiable functions from  $\RR^{2d}$ into $\RR^d$ and $\RR^{d\times r}$, for some $r\geq 1$; and
 let  $W_t$  be an $r$-dimensional Brownian motion.
   For any $\mu\in \Pa_2(\RR^d)$ and any time horizon $s\geq 0$ we denote by
  $X_{s,t}^{\mu}(x)$ be the stochastic flow defined for any $t\in [s,\infty[$ and any starting point $x\in \RR^d$ by the McKean-Vlasov diffusion 
    \begin{equation}\label{diff-st-ref-general}
    dX^{\mu}_{s,t}(x)=b_t\left(\phi_{s,t}(\mu),X^{\mu}_{s,t}(x)\right)~dt+\sigma_{t}\left(\phi_{s,t}(\mu),X^{\mu}_{s,t}(x)\right)~
    dW_t
    \end{equation}
In the above display, $\phi_{s,t}$ stands for the evolution semigroup 
    $$
     \phi_{s,t}(\mu)(dy)=\mu P^{\mu}_{s,t}(dy):=\int~\mu(dx)~P^{\mu}_{s,t}(x,dy)\quad\mbox{\rm with}\quad P^{\mu}_{s,t}(x,dy):=\PP(X^{\mu}_{s,t}(x)\in dy)
    $$
We further assume that the mean field drift and diffusion functions are given by 
    $$
    b_t\left(\eta,y\right):=\int~\eta(dx)~b_t(x,y)
    \quad\mbox{\rm and}\quad
    \sigma_{t}\left(\eta,y\right):=\int~\eta(dx)~\sigma_{t}(x,y)
    $$ 

    We shall assume that the nonlinear diffusion flow (\ref{diff-st-ref-general}) is well defined.
For instance, the existence of this flow is ensured as son as
$b_t$ and $\sigma_t$ are Lipschitz, see for instance~\cite{graham,huang}. 

    The mean field particle system associated with (\ref{diff-st-ref-general})  is defined 
     by the stochastic flow  $\xi_{s,t}(z)=(\xi_{s,t}^i(z))_{1\leq i\leq N}$ of a system of  $N$ interacting diffusions
 \begin{equation}\label{diff-st-ref-general-mf}
 d\xi_{s,t}^i(z)=b_t\left(m(\xi_{s,t}^i(z)),\xi_{s,t}^i(z)\right)~dt+  \sigma_{t}\left(m(\xi_{s,t}^i(z)),\xi_{s,t}^i(z)\right)~ dW^{i}_t 
     \end{equation}
 with the empirical measures
 $$
 m(\xi_{s,t}^i(z)):=\frac{1}{N}\sum_{1\leq j\leq N}\delta_{\xi_{s,t}^j(z)} 
 $$
 In the above displayed formulae,  $\xi_{s,s}(z)=z=(z^i)_{1\leq N}\in (\RR^d)^N$ stands for the initial configuration and $W^{i}_t$
 are $N$ independent copies of $W_t$.

\subsection{Statement of some main results and article organisation}

 To motivate this study, the variational calculus developed in the article is illustrated with the following example
  \begin{equation}\label{ex-langevin}
r=d\qquad \sigma(x,y)=\sigma_0~I\quad \mbox{\rm and}\quad b(x,y)=-\nabla U(y)-\nabla V(y-x)
 \end{equation}
  for some $\sigma_0>0$, some confinement type potential function $U$ (a.k.a. the exterior potential) and some interaction potential function $V$.
  This class of nonlinear diffusions and the corresponding particle interpretations were introduced by  H. P. McKean in~\cite{mckean-1,mckean-2}.
  The extended versions of these Langevin-type nonlinear diffusions on Riemannian manifolds are discussed in the end of section~\ref{sec-manifold} as well as in section~\ref{subsec4.3}.

 Nonlinear diffusions (\ref{diff-st-ref-general}) with constant diffusion and gradient-type drifts  (\ref{ex-langevin}) arise in fluid mechanics, and more particularly
  in the modeling of granular flows~\cite{bene-1,bene-2,toscani,villani}.  In this context, $ \phi_{s,t}$ represents the evolution semigroup of the 
  velocity of a diffusive particule interacting with the distribution of the particles around its location and following some confinement exterior potential.
 In this interpretation, the mean field particle model (\ref{diff-st-ref-general-mf}) can be seen as a particle-type representation of the granular flow.

  In the last two decades, the analysis of the long time behavior of this particular class of gradient type flow diffusions have been developed in various directions:

The first articles on the long time behavior of these models are the couple of articles by Tamura~\cite{tamura,tamura-2}. 
The stability properties of one dimensional models has been started in~\cite{benachour-1,benachour-2} as well as in~\cite{bene-1}, see also \cite{calvez,c-toscani,toscani}. 

Since this period, several sophisticated probabilistic techniques have been developed to analyze the long time behavior of these Langevin-type nonlinear diffusions, including log-Sobolev functional inequalities~\cite{malrieu,malrieu-2}, entropy dissipation~\cite{cmcv-2006,cordero},
as well as gradient flows in Wasserstein metric spaces and optimal transportation inequalities~\cite{gbg-13,cmcv-2006,cattiaux,otto}, combining
the functional $\Gamma_2$  Bakry-Emery method~\cite{bakry}, with the Otto-Villani approach~\cite{otto-2}.  The long time self-stabilizing behavior of  this class of processes in multi-wells landscapes has  also been developed by J. Tugaut in a series of articles~\cite{tugaut-1,tugaut-2,tugaut-4,tugaut-5,tugaut-6}.
For a more thorough discussion on this subject we refer to the recent article~\cite{tugaut-hong}, and the references therein. 

Unfortunately, most of the probabilistic techniques discussed above only apply to gradient flow type diffusions of the form (\ref{ex-langevin}). 
 The variational calculus developed in the present article is not
 restricted to this class  of gradient-type nonlinear models. 
 Nevertheless, because of their importance in practice this introduction illustrates some of our main results in this context.
 
 Firstly, and rather surprisingly, the variational methodology developed in the present article applies directly
 to gradient flow models of the form  (\ref{ex-langevin}), simplifying considerably both of
their stability analysis as well as the convergence analysis of their mean field particle interpretations.  

This framework
also allows to relax unnecessary technical conditions such as the symmetry of the interaction potential function, or the
invariance of the center of mass,  currently used in the literature on this subject (see for instance~\cite{tamura}, as well as section 2 in~\cite{cmcv-2006}, and section 1 in the more recent article~\cite{gbg-13}).  It also allows to derive uniform as well as almost sure exponential stability inequalities at the level of the nonlinear diffusion flow. 
For instance, 
when $V$ is an even convex function with bounded Hessian $\Vert\nabla^2V\Vert_2:=\sup_x\Vert\nabla^2V(x)\Vert_2<\infty$, and when $\nabla^2 U\geq \lambda~I$, for some $\lambda>0$ we have the almost sure estimates
  \begin{equation}\label{ref-intro-illustration}
\Vert X^{\eta}_{s,t}(x)-X^{\mu}_{s,t}(y)\Vert\leq \Vert\nabla^2V\Vert_2~(t-s)~e^{-\lambda(t-s)}~\WW_2(\eta,\mu)+e^{-\lambda(t-s)}~~\Vert x-y\Vert
 \end{equation}
The above estimate is also met for odd interaction potential, as soon as $\nabla^2 U(y)+\nabla^2 V(y-x)\geq \lambda~I$. In the above display, it is implicitly assumed that the stochastic flows are driven by the same Brownian motion.

These almost sure inequalities are direct consequence of the contraction inequality (\ref{ref-nablax-estimate-0-ae-bis}), the remark
(\ref{ref-ex-langevin}) and the almost sure estimates stated in corollary~\ref{ref-intro-theo}.

To the best of our knowledge, the almost sure exponential decays (\ref{ref-intro-illustration}) are the first result of this type for this  class of nonlinear gradient flow diffusions. 

Consider a pair of random variables $(Z_0,Z_1)$ with distributions $(\mu_0,\mu_1)$ on $\RR^d$ and set 
 \begin{equation}\label{def-mu-epsilon}
Z_{\epsilon}:=
(1-\epsilon)~Z_0+\epsilon~Z_1\qquad\mu_{\epsilon}:=\mbox{\rm Law}(Z_{\epsilon})\quad \mbox{\rm and}\quad 
X^{\epsilon}_{s,t}:=X^{\mu_{\epsilon}}_{s,t}(Z_{\epsilon})
 \end{equation}
Under the assumptions on the potential functions discussed above, for any differentiable function $f$ on $\RR^d$ with bounded gradient we have the first order differential formula
 \begin{equation}\label{diff-epsilon-ref-intro}
 \left[\phi_{s,t}(\mu_1)-   \phi_{s,t}(\mu_0)\right](f)=\int_0^1~ \partial_{\epsilon} \phi_{s,t}(\mu_{\epsilon})(f)~d\epsilon
 \end{equation}
 with the linear differential operator
$$
  \partial_{\epsilon} \phi_{s,t}(\mu_{\epsilon})(f):=
  \EE\left(\left\langle\partial_{\epsilon}X^{\epsilon}_{s,t},\nabla f(X^{\epsilon}_{s,t})\right\rangle\right)\quad \mbox{\rm s.t.}\quad
  \vert  \partial_{\epsilon} \phi_{s,t}(\mu_{\epsilon})(f)\vert\leq e^{-\lambda (t-s)}~\Vert\nabla f\Vert
  $$
For a more precise statement we refer to theorem~\ref{stab-nl-theo}. Almost sure and uniform estimates of the first order differential
maps $\epsilon\mapsto \partial_{\epsilon}X^{\epsilon}_{s,t}$ are also provided in theorem~\ref{theo-as-ref}.

Section~\ref{stab-particle-sec} also presents a differential calculus to estimate the gradient $\nabla\xi_{s,t}(z)$ of the stochastic flow $\xi_{s,t}(z)$ of the interacting particle model (\ref{diff-st-ref-general-mf}). Under the assumptions on the potential functions discussed above, we shall prove the following uniform spectral norm estimate
$$
\Vert\nabla\xi_{s,t}(z)\Vert_2\leq e^{-\lambda (t-s)}
$$
The above result is a direct consequence of theorem~\ref{theo-Aa}. The above estimate ensures that the $N$-particle model converges exponentially fast to its invariant measure with some exponential decay that doesn't depends on the number of particles. The latter property can also be checked
using more sophisticated Logarithmic Sobolev inequalities~\cite{malrieu}. To the best of our knowledge, the almost sure exponential decays stated above are the first result of this type for this class of interacting Langevin-type diffusions.    

Section~\ref{prop-chaos-sec} also provides a natural differential calculus to derive quantitative and uniform propagation of chaos estimates for nonlinear diffusions of the form (\ref{diff-st-ref-general-mf}). Applying these results to interacting Langevin-type diffusions, without further work we recover the uniform estimates stated in theorem 1.2 in~\cite{malrieu}.

We emphasize that the differential calculus presented in this article allows to consider nonlinear diffusions evolving in differential manifolds.
 This should not come as a surprise since our framework allows to enter the variations of
the diffusion matrices associated with these stochastic models which encapsulates the Riemannian structure of the manifold. 

We illustrate these comments in the end of section~\ref{sec-1st-o-d} with a rather detailed discussion of an elementary nonlinear geometric-type diffusion.  The manifold version of (\ref{diff-epsilon-ref-intro}) is also provided in theorem~\ref{diff-contraction-general-g}.

We also underline that the variational calculus on differentiable manifolds developed in section~\ref{sec-manifold} provides another view and additional results for the diffusions in $\RR^d$ endowed, when possible, with the Riemannian metric under which these diffusions are Brownian motion with drift. In this context,  different types of synchronous coupling  lead to gradient flow estimates where gradients of the diffusion functions are replaced by Ricci curvatures. 

Quantitative propagation of chaos estimates of mean field particle systems on Riemannian manifolds are provided in section~\ref{subsec4.3}. Special attention is paid to derive uniform estimates w.r.t. the time horizon.

 \section{Nonlinear diffusion semigroups}
 \subsection{Some gradient flow estimates}   \label{sec-gfe}
This section presents some basic properties of the first variational equation associated with the nonlinear diffusion (\ref{diff-st-ref-general}).
   Let $\sigma_{k,t}$ be the $k$-th column vector of $\sigma_{t}$, 
and let   $\nabla_{u}b_t(x,y)$ and $\nabla_{u}\sigma_{k,t}(x,y)$ be the  gradient of the functions $b_t(x,y)$ and $\sigma_{k,t}(x,y)$
w.r.t.
 the coordinate $u\in\{x,y\}$. We also let $X^{i,\mu}_{s,t}(x)$ be the $i$-th coordinate of the column vector $X^{\mu}_{s,t}(x)$. 
        The Jacobian $\nabla X^{\mu}_{s,t}(x)$ of the diffusion flow $X^{\mu}_{s,t}(x)$ is given by the gradient $(d\times d)$-matrix
$$
\begin{array}{l}
\displaystyle\nabla X^{\mu}_{s,t}(x):=\left(\nabla X^{1,\mu}_{s,t}(x),\ldots,\nabla X^{d,\mu}_{s,t}(x)\right)
\\
\\
\displaystyle\Longrightarrow
d \,\nabla X^{\mu}_{s,t}(x)=\nabla X^{\mu}_{s,t}(x)~\left[\nabla b_t\left(\phi_{s,t}(\mu),X^{\mu}_{s,t}(x)\right)~dt+\sum_{1\leq k\leq r}\nabla \sigma_{t,k}\left(\phi_{s,t}(\mu),X^{\mu}_{s,t}(x)\right)~dW^k_t\right]
\end{array}
$$
 
Consider the regularity condition stated below:\\
 
 {\em $(H_A)$ : There exists some $\lambda_A\in \RR$ such that for any $x,y\in \RR^d$ and $t\geq 0$ we have 
   \begin{equation}\label{ref-mat-Upsilon}
\qquad A_t(x,y):=\nabla_y b_t(x,y)+\nabla_y b_t(x,y)^{\prime}+\sum_{1\leq k \leq r}\nabla_y \sigma_{k,t}(x,y)\nabla_y \sigma_{k,t}(x,y)^{\prime}\leq -2\lambda_A~I
 \end{equation} 
 }

 This spectral condition produces several gradient estimates. For instance, we have the following uniform estimate
  \begin{equation}\label{ref-nablax-estimate-0}
(H_A)\Longrightarrow \EE\left(\Vert\nabla X^{\mu}_{s,t}(x)\Vert_2^2\right)^{1/2}\leq \EE\left(\Vert\nabla X^{\mu}_{s,t}(x)\Vert_F^2\right)^{1/2}\leq \sqrt{d}~e^{-\lambda_A(t-s)}
 \end{equation}
 In addition, we have the almost sure estimate
  \begin{equation}\label{ref-nablax-estimate-0-ae}
(H_A)\quad \mbox{\rm and}\quad \nabla_y \sigma_{k,t}(x,y)=0\Longrightarrow \Vert\nabla X^{\mu}_{s,t}(x)\Vert_2\leq e^{-\lambda_A(t-s)}
 \end{equation}
The proofs of the above assertions are provided in the appendix, on page~\pageref{ref-nablax-estimate-proof}. 
 For the nonlinear Langevin diffusion discussed in (\ref{ex-langevin}) we have
   \begin{equation}\label{ref-HA-Langevin}
(H_A)\quad \Longleftrightarrow\quad  \nabla^2 U(y)+\nabla^2 V(y-x)\geq \lambda_A~I\quad \Longrightarrow  \Vert\nabla X^{\mu}_{s,t}(x)\Vert_2\leq e^{-\lambda_A(t-s)}
 \end{equation}

Arguing as in (\ref{taylor-first}) we readily check the following proposition.
\begin{prop}
Assume $(H_A)$ is satisfied. In this situation, we have 
  \begin{equation}\label{ref-nablax-estimate-0-bis}
 \EE\left(\Vert X^{\mu}_{t}(x)-X^{\mu}_{t}(y)\Vert^2\right)^{1/2}\leq \sqrt{d}~e^{-\lambda_A(t-s)}~~\Vert x-y\Vert
 \end{equation}
In addition, we have the almost sure estimate
  \begin{equation}\label{ref-nablax-estimate-0-ae-bis}
 \nabla_y \sigma_{k,t}=0\Longrightarrow \Vert X^{\mu}_{t}(x)-X^{\mu}_{t}(y)\Vert\leq e^{-\lambda_A(t-s)}~~\Vert x-y\Vert
 \end{equation}
 \end{prop}
 Whenever $\lambda_A<0$  the above estimates ensure that the transition semigroup $P^{\mu}_{s,t}$ is exponentially stable, that is we have that
   \begin{equation}\label{ref-eta-mu-cv}
   \WW_2\left(\eta_0 P^{\mu}_{s,t},\eta_1 P^{\mu}_{s,t}\right)\leq c~\exp{\left[-\lambda_A(t-s)\right]}~  \WW_2\left(\eta_0 ,\eta_1 \right)
 \end{equation}
  These contraction inequalities quantify the stability of the stochastic flow $X^{\mu}_{s,t}(x)$  w.r.t. the initial state $x$, but
they don't give any information of the stability properties of the nonlinear semigroup $\phi_{s,t}(\mu)$ w.r.t. the initial measure $\mu$.
 \subsection{A first order differential calculus}\label{sec-1st-o-d}
This section presents a natural  first order differential calculus to analyze the  
stability properties of the nonlinear semigroup $\phi_{s,t}(\mu)$. Consider  the matrices
\begin{equation}\label{def-B-Sigma}
B_{t}(z_1,z_2):=\left[
\begin{array}{cc}
\nabla_yb_t\left(z_2,z_1\right)&\nabla_xb_t\left(z_1,z_2\right)\\
&\\
\nabla_xb_t\left(z_2,z_1\right)&\nabla_yb_t\left(z_1,z_2\right)
\end{array}
\right]
\quad
D_t:=\sum_{1\leq k\leq r}\left[
\begin{array}{cc}
\nabla_x \sigma_{t,k}~\nabla_x \sigma_{t,k}^{\prime}&\nabla_x \sigma_{t,k}~\nabla_y \sigma_{t,k}^{\prime}\\
&\\
\nabla_y \sigma_{t,k}~\nabla_x \sigma_{t,k}^{\prime}&\nabla_y \sigma_{t,k}~\nabla_y \sigma_{t,k}^{\prime}
\end{array}
\right]
 \end{equation}
 
 In this notation, our second regularity condition takes the following form:\\

 {\em $(H_C)$ : There exists some $\lambda_C\in \RR$ such that for any $x,y\in \RR^d$ and $t\geq 0$ we have 
 \begin{equation}\label{hyp-b-sigma-general}
C_t(x,y):= \frac{1}{2}\left[B_{t}\left(x,y\right)+ B_{t}\left(x,y\right)^{\prime}\right]+
D_t\left(x,y\right)\leq -\lambda_C~I
 \end{equation}}
 
Let $Z_{\epsilon}$ be the collection of random variables with distribution $\mu_{\epsilon}$ defined in (\ref{def-mu-epsilon}).
We also consider a couple of independent stochastic flows
\begin{equation}\label{def-Yepsilon}
X^{\epsilon}_{s,t}:=X^{\mu_{\epsilon}}_{s,t}(Z_{\epsilon})\quad \mbox{\rm and}\quad
Y^{\epsilon}_{s,t}:=Y^{\mu_{\epsilon}}_{s,t}(\overline{Z}_{\epsilon})
\end{equation}
driven by independent 
Brownian motions, say $W_t=(W^k_t)_{1\leq k\leq d}$ and $\overline{W}_t=(\overline{W}^k_t)_{1\leq k\leq d}$, and starting from a
couple of independent random variables $Z_{\epsilon}$ and $\overline{Z}_{\epsilon}$ with the same law.

In the further development of this section, we denote by ${\EE}_{X}(\point)$ the expectation operator w.r.t. the Brownian motion
$W_t=({W}^k_t)_{1\leq k\leq d}$ and the random variable  $Z_{\epsilon}$. In this notation, we have
$$
 dY^{\epsilon}_{s,t}=\EE_X\left[b_t\left(X^{\epsilon}_{s,t},Y^{\epsilon}_{s,t}\right)\right]~dt+\EE_X\left[
 \sigma_{t}\left(X^{\epsilon}_{s,t},Y^{\epsilon}_{s,t}\right)\right]~d\overline{W}_t
$$
This implies that
  \begin{equation}\label{def-partial-Y}
\begin{array}{l}
 \displaystyle d\left[\partial_{\epsilon}Y^{\epsilon}_{s,t}\right]
=\EE_X\left[\nabla_xb_t\left(X^{\epsilon}_{s,t},Y^{\epsilon}_{s,t}\right)^{\prime}\partial_{\epsilon}X^{\epsilon}_{s,t}
 +\nabla_yb_t\left(X^{\epsilon}_{s,t},Y^{\epsilon}_{s,t}\right)^{\prime}\partial_{\epsilon}Y^{\epsilon}_{s,t}\right]~dt\\
 \\
 \hskip3cm \displaystyle+\sum_{1\leq k\leq r}\EE_X\left[
\nabla_x \sigma_{t,k}\left(X^{\epsilon}_{s,t},Y^{\epsilon}_{s,t}\right)^{\prime}\partial_{\epsilon}X^{\epsilon}_{s,t}+\nabla_y \sigma_{t,k}\left(X^{\epsilon}_{s,t},Y^{\epsilon}_{s,t}\right)^{\prime}\partial_{\epsilon}Y^{\epsilon}_{s,t}\
\right]~d\overline{W}^k_t
\end{array}
 \end{equation}
 with the initial condition $$\partial_{\epsilon}Y^{\epsilon}_{s,s}=\partial_{\epsilon}\overline{Z}_{\epsilon}=\overline{Z}_1-\overline{Z}_0$$
 A simple calculation yields the following estimate
  \begin{equation}\label{diff-contraction-general-pre}
 \displaystyle \partial_t\,\EE\left[\left\Vert\partial_{\epsilon}Y^{\epsilon}_{s,t}\right\Vert^2\right]\leq \EE\left(
\left[\partial_{\epsilon}X^{\epsilon}_{s,t},\partial_{\epsilon}Y^{\epsilon}_{s,t}\right]^{\prime}C_{t}\left(X^{\epsilon}_{s,t},Y^{\epsilon}_{s,t}\right)
\left[
\begin{array}{c}
\partial_{\epsilon}X^{\epsilon}_{s,t}\\
\partial_{\epsilon}Y^{\epsilon}_{s,t}
\end{array}
\right]\right) \end{equation}
The inequality in the above display can be turned into an equality when $D_t=0$.
Also note that
   \begin{eqnarray*}
(H_C)\Longrightarrow \EE\left(\Vert Y^{0}_{s,t}-Y^{1}_{s,t}\Vert^2\right)&\leq& \int_0^{1} \EE\left(\left\Vert\partial_{\epsilon}Y^{\epsilon}_{s,t}\right\Vert^2\right)~d\epsilon\leq e^{-2\lambda_C (t-s)}~\EE\left(\Vert \overline{Z}_1-\overline{Z}_0\Vert^2\right)
    \end{eqnarray*}
    Let  $ \Ca^1_b(\RR^d)$ be the set of differentiable functions on $\RR^d$ with bounded derivative. 
 A direct consequence of the fundamental theorem of calculus yields the following theorem.
\begin{theo}\label{stab-nl-theo}
For any $s\leq t$ and any $f\in \Ca^1_b(\RR^d)$ and  $\mu_0,\mu_1\in \Pa_2(\RR^d)$  we have the first order differential formula
 (\ref{diff-epsilon-ref-intro}).
 In addition, we have the exponential contraction inequality
           \begin{equation}\label{diff-contraction-general}
(H_C)\quad\Longrightarrow\quad \WW_2\left(\phi_{s,t}(\mu_0),\phi_{s,t}(\mu_1)\right)\leq e^{-\lambda_C (t-s)}~\WW_2(\mu_0,\mu_1)
 \end{equation}

   \end{theo}

When $\lambda_C>0$, the above theorem provides an 
alternative and rather elementary proof of the exponential asymptotic stability of time varying McKean-Vlasov diffusions with non necessarily homogenous diffusion functions. To the best of our knowledge this stability property is the first result of this type for this general class of nonlinear diffusions.
 
 For the Langevin-type diffusion discussed in (\ref{ex-langevin}) we have $D_t=0$ and the matrix $C_t$ reduces to
$$
\begin{array}{l}
 -C_{t}(z_1,z_2)
 =\left[
\begin{array}{cc}
\nabla^2 U(z_1)&0\\
&\\
0&\nabla^2 U(z_2)
\end{array}
\right]\\
\\
\hskip3cm+\left[
\begin{array}{cc}
\nabla^2 V(z_1-z_2)&\displaystyle-\frac{\left[\nabla^2 V(z_2-z_1)+\nabla^2 V(z_1-z_2)\right]}{2}\\
&\\
\displaystyle-\frac{\left[\nabla^2 V(z_2-z_1)+\nabla^2 V(z_1-z_2)\right]}{2}&\nabla^2 V(z_2-z_1)
\end{array}
\right]
\end{array}
$$
When $V$ is odd we have 
\begin{equation}\label{ref-ex-langevin-odd}
 (H_C)~\Longleftrightarrow~ \nabla^2 U(z_1)+\nabla^2 V(z_1-z_2)\geq \lambda_C~I~\Longleftrightarrow~(H_A)
\end{equation}
In the reverse angle, if $V$ is even and convex then we have
\begin{equation}\label{ref-ex-langevin}
 (H_C)~\Longleftrightarrow~ \nabla^2 U\geq \lambda_C~I\Longrightarrow (H_A)
\end{equation}
As expected, explicit formulae are available for linear and Gaussian models. For instance, when  
$$
 b_t(x,y)=A_1\,x+A_2\,y\quad \mbox{\rm and}\quad \sigma_t(x,y)=R^{1/2}\quad \mbox{\rm with}\quad A_1,A_2\in\RR^{d\times d}\quad \mbox{\rm and}\quad 
R\geq 0   
$$
 the diffusion flow  $X_{s,t}^{\mu}(x)\in \RR^d$ 
is linear w.r.t.  $\mu$ and given for any $x\in \RR^d$ by the formula
$$
 X^{\mu}_{s,t}(x)=e^{A_2(t-s)}(x-\mu(e))+e^{[A_1+A_2](t-s)}~\mu(e)+\int_s^te^{A_2(t-u)}R^{1/2}~dW_u
  $$
 In the above display,  $e(x)=x$ stands for the identity function on $\RR^d$.  In this context, the process $X^{\epsilon}_{s,t}$ defined in (\ref{def-Yepsilon}) is also given by the formula
$$
 \begin{array}{l}
\displaystyle
 X^{\epsilon}_{s,t}=e^{A_2(t-s)}(Z_{\epsilon}-\mu_{\epsilon}(e))+e^{[A_1+A_2](t-s)}~\mu_{\epsilon}(e)+\int_s^te^{A_2(t-u)}R^{1/2}~dW_u\\
 \\
  \displaystyle\Longrightarrow\partial_{\epsilon}X^{\epsilon}_{s,t}=e^{A_2(t-s)}((Z_1-Z_0)-\EE(Z_1-Z_0))+e^{[A_1+A_2](t-s)}~\EE(Z_1-Z_0)
 \end{array}$$ 
 This yields the rather crude estimate
 $$
 \EE\left(\Vert \partial_{\epsilon}X^{\epsilon}_{s,t} \Vert^2\right)\leq \left[\Vert e^{[A_1+A_2](t-s)}\Vert_2^2+\Vert  e^{A_2(t-s)}\Vert_2^2\right]~\EE(\Vert Z_1-Z_0\Vert^2)
 $$ Up to some constant, this shows that  the r.h.s. Wasserstein contraction estimate in (\ref{diff-contraction-general}) is met with $-\lambda_C=\rho(A_1+A_2)\vee \rho(A_2)$. Applying Coppel's inequality (cf. Proposition 3 in~\cite{coppel}) we can also choose
 $-\lambda_C=\left[\varsigma(A_1+A_2)\vee \varsigma(A_2)\right](1-\delta)$ for any $0<\delta<1$,  where $\varsigma(A):=\max_i{\left\{\mbox{\rm Re}[\lambda_i(A)]\right\}}\leq \rho(A)$ stands for the spectral abscissa of a square matrix $A$. 
 
It may happen the stochastic flow (\ref{diff-st-ref-general}) remains in some domain $S\subset \RR^d$.
The simplest model we have in head is the 
 geometric diffusion on $S=[0,\infty [$ associated with the parameters
$$
   b_t(x,y)=[a_1-a_2~x]~y\quad \mbox{\rm and}\quad \sigma_t(x,y)=\sigma_0~y\quad \mbox{\rm with}\quad a_1\in\RR\quad \mbox{\rm and}\quad 
   a_2,\sigma_0> 0   
$$
 In this situation, the diffusion flow  $X_{s,t}^{\mu}(x)\in S$ 
is nonlinear w.r.t.  $\mu$ and given for any $x\in S$  by
\begin{equation}\label{geo-diff}
X_{s,t}^{\mu}(x)=\psi_{t-s}(\mu)~\Ea_{s,t}(W)~x\quad\mbox{\rm with}\quad \Ea_{s,t}(W):=\exp{\left[\sigma_0(W_t-W_s)-\frac{\sigma_0^2}{2}(t-s)\right]}
\end{equation}
with the function $\psi_{t}$ defined by
$$
\psi_{t}(\mu)=
\displaystyle\frac{1}{e^{-a_1t}+a_2\,\mu(e)~\theta_{a_1}(t)}\quad \mbox{\rm with}\quad \theta_{a_1}(t):=a_1^{-1}(1-e^{-a_1t})
$$
 In the above display, we have used the convention $\theta_{0}(t)=t$.  In this context, the process $X^{\epsilon}_{s,t}$ defined in (\ref{def-Yepsilon}) is also given by the formula
 $$
 \begin{array}{l}
\displaystyle
 X^{\epsilon}_{s,t}=\psi_{t-s}(\mu_{\epsilon})~\Ea_{s,t}(W)~Z_{\epsilon}\\
 \\
 \displaystyle\Longrightarrow\partial_{\epsilon}X^{\epsilon}_{s,t}=\psi_{t-s}(\mu_{\epsilon})~\Ea_{s,t}(W))~\left[(Z_1-Z_0)-a_2\,\theta_{a_1}(t)~\psi_{t-s}(\mu_{\epsilon})~Z_{\epsilon}~\EE(Z_1-Z_0)\right]
 \end{array}$$
Assume that $a_1<0$ is chosen so that $\vert a _1\vert> \sigma_0^2/2$. In this situation, for any  $x,y\in S$ we have
  $$
 A_t(x,y)=2[a_1-a_2~x]+\sigma_0^2\leq 2a_1+\sigma_0^2\Longrightarrow (H_A)\quad \mbox{\rm with}\quad\lambda_{A}=\vert a_1\vert-\sigma_0^2/2<0
$$
as well as
 $$
\psi_{t}(\mu)=\frac{\vert a_1\vert~e^{-\vert a_1\vert t}}{\vert a_1\vert+a_2\,\mu(e)~\left(1-e^{-\vert a_1\vert t}\right)} \leq e^{-\vert a_1\vert t}
 $$
This yields the estimate
$$
 \EE\left[ [\partial_{\epsilon}X^{\epsilon}_{s,t}]^2\right]\leq \left[1+\vert a^{-1}_1a_2\vert~e^{-\vert a_1\vert (t-s)}~\left(\EE(Z_0^2)^{1/2}\vee \EE(Z_1^2)^{1/2}\right)\right]^2~e^{-(2\vert a_1\vert-\sigma_0^2)(t-s)}~ \EE((Z_1-Z_0)^2)
$$
Up to some constant, this shows that  the r.h.s. Wasserstein contraction estimate in (\ref{diff-contraction-general}) is met with $\lambda_C=\vert a_1\vert-\sigma_0^2/2$.

The analysis of nonlinear diffusions on more general differentiable manifolds is based on more sophisticated differential techniques. 
The extension of the variational calculus developed above to this class of stochastic processes on manifolds is provided in
section~\ref{sec-manifold}.

We end this section with some illustrations of our results on time homogeneous models $(b_t,\sigma_t)=(b,\sigma)$ satisfying condition $(H_C)$. We
set
$\phi_t:=\phi_{0,t}$, and $P^{\mu}_t:=P^{\mu}_{0,t}
$. 
By theorem~\ref{stab-nl-theo},   there exists an unique invariant measure $$\pi=\phi_{t}(\pi)\quad\mbox{\rm and}\quad \WW_2\left(\phi_t(\mu),\pi\right)\leq e^{-\lambda_Ct}~  \WW_2\left(\mu,\pi \right)
   $$
For the nonlinear Langevin diffusion discussed in (\ref{ex-langevin}) condition $(H_C)$ is met when (\ref{ref-ex-langevin-odd}) or (\ref{ref-ex-langevin}) are satisfied. In this context, $X^{\pi}_t:=X^{\pi}_{0,t}$ is a conventional Langevin diffusion given by the time homogeneous  stochastic differential equation
$$
   dX^{\pi}_t=-\nabla V^{\pi}(X^{\pi}_t)~dt+\sigma_0~
    dW_t\quad \mbox{\rm with}\quad  2^{-1}V^{\pi}(y)=U(y)+\int\pi(dx)~V(y-x)
$$
In this situation, the unique invariant measure of $X^{\pi}_t$ is given by 
$$
\varpi(\pi)(dx):=\frac{1}{v_{\pi}}\exp{\left[-\frac{1}{\sigma_0}~V^{\pi}(x)\right]}~dx
\quad \mbox{\rm with}\quad  v_{\pi}:=\int~\exp{\left[-\frac{1}{\sigma_0}~V^{\pi}(x)\right]}~dx
$$
In the above display, $dx$ stands for the Lebesgue measure on $\RR^d$. In this case the measure $\pi=\phi_{t}(\pi)=\pi P^{\pi}_t$ is the unique solution of the equation $\pi=\varpi(\pi)$.
We underline that the uniqueness of the invariant measure is not ensured for double-well confinement potential functions and too small noise. Further details on this subject including a description of the invariant measures for small noise can be found in the series of articles~\cite{h-t-2010,2-h-t-2010,3-h-t-2010}.
 
Whenever $(H_C)$ is met,  we also have the uniform moment estimates
\begin{equation}
 \label{unif-moments}
\phi_t(\mu)(\Vert e\Vert^2)^{1/2}\leq \pi(\Vert e\Vert^2)^{1/2}+\WW_2\left(\mu,\pi \right)
\end{equation}
In the same vein, when when $(H_A)$ and $(H_C)$ are met  we have
$$
\EE\left[\Vert X^{\mu}_t(x)\Vert^{2}\right]^{1/2}\leq \pi(\Vert e\Vert^2)^{1/2}+ \WW_2\left(\delta_x P^{\mu}_{t},\pi\right)\leq 
c~\left[\pi(\Vert e\Vert^2)\vee \mu(\Vert e\Vert^2)\right]^{1/2}~\left[1+\Vert x\Vert\right]
$$
for some finite constant $c$.
The last assertion comes from the fact that
\begin{eqnarray*}
 \WW_2\left(\delta_x P^{\mu}_{t},\pi\right)&\leq &  \WW_2\left(\delta_x P^{\mu}_{t},\phi_{t}(\mu)\right)+ \WW_2\left(\phi_{t}(\mu),\pi\right) \leq c~e^{-(\lambda_A\wedge\lambda_C)t}~  \left[\WW_2\left(\delta_x ,\mu \right)+ \WW_2\left(\mu,\pi \right)\right]
\end{eqnarray*}

 \subsection{Some almost sure estimates}
  We fix the parameters $\epsilon$ and some given time horizon $s\geq 0$, and we set $y_{t}:=\partial_{\epsilon}Y^{\epsilon}_{s,t}$, for any $t\in [s,\infty[$,     with the process  $Y^{\epsilon}_{s,t}$ defined in  (\ref{def-partial-Y}).  Also consider the processes
$$
dz_{t}:=z_{0,t}~dt+\sum_{1\leq k\leq r}~z_{k,t}~d\overline{W}^k_t\quad \mbox{\rm and}\quad
d\Za_{t}:=   \Za_{0,t}~dt+\sum_{1\leq k\leq r}   \Za_{k,t}~d\overline{W}^k_t
$$
with  the collection of processes
   \begin{eqnarray*}
z_{0,t}&:=&\EE_X\left[\nabla_xb_t\left(X^{\epsilon}_{s,t},Y^{\epsilon}_{s,t}\right)^{\prime}\partial_{\epsilon}X^{\epsilon}_{s,t}\right]~
\qquad
z_{k,t}:=\EE_X\left[
\nabla_x \sigma_{t,k}\left(X^{\epsilon}_{s,t},Y^{\epsilon}_{s,t}\right)^{\prime}\partial_{\epsilon}X^{\epsilon}_{s,t}
\right]
\\
   \Za_{0,t}&:=&\EE_X\left[\nabla_yb_t\left(X^{\epsilon}_{s,t},Y^{\epsilon}_{s,t}\right)^{\prime}\right]\qquad \mbox{\rm and}\qquad
    \Za_{k,t}:=\EE_X\left[\nabla_y \sigma_{t,k}\left(X^{\epsilon}_{s,t},Y^{\epsilon}_{s,t}\right)^{\prime}\right]
    \end{eqnarray*}
 In this notation,
the evolution equation (\ref{def-partial-Y}) reduces to
$$
dy_{t}:=dz_{t}+d\Za_{t}~y_{t}
$$
Let
 $t\in [s,\infty[\mapsto
 \Ea_{t}
 $ be the solution of the matrix evolution equation
 $$
 \begin{array}{l}
d\Ea_{t}:=d\Za_{t}~\Ea_{t}\quad \mbox{\rm with}\quad \Ea_{s}=I\quad \mbox{\rm and set}\quad
 \forall t\in [u,\infty[\qquad \Ea_{u,t}:=\Ea_{t}\Ea_{u}^{-1}\\
\\
\Longrightarrow \forall t\in [u,\infty[\qquad  d\, \Ea_{u,t}:=d\Za_{t}~\Ea_{u,t}
\end{array} $$
  In this notation, we readily check that
 $$
 y_{t}=\Ea_{s,t}~ y_s+\int_s^t\Ea_{u,t}\left(dz_{u}-\sum_{1\leq k\leq r}   \Za_{k,u}z_{k,u}~du\right)
 $$
Whenever condition $(H_A)$ is met,  for any given $u\geq 0$ and any $t\in [u,\infty[$ we have
\begin{eqnarray*}
\displaystyle  d \left[\Ea_{u,t}^{\prime}\, \Ea_{u,t}\right]&=&\Ea_{u,t}^{\prime}\left[  \Za_{0,t}+\Za_{0,t}^{\prime}+\sum_{1\leq k\leq r}   \Za_{k,t}^{\prime} ~\Za_{k,t}\right]\Ea_{u,t}~dt+\sum_{1\leq k\leq r}\Ea_{u,t}^{\prime}\left(   \Za_{k,t}+   \Za_{k,t}^{\prime}\right)\Ea_{u,t}~d\overline{W}^k_t\\
&\leq &-2\lambda_A~\Ea_{u,t}^{\prime}\, \Ea_{u,t}~dt+\sum_{1\leq k\leq r}\Ea_{u,t}^{\prime}\left(   \Za_{k,t}+   \Za_{k,t}^{\prime}\right)\Ea_{u,t}~d\overline{W}^k_t
\end{eqnarray*}  
  This shows that
  $$
  (H_A)\quad \mbox{\rm and}\quad \nabla_y \sigma_{k,t}=0\Longrightarrow \Ea_{u,t}^{\prime}\, \Ea_{u,t}\leq e^{-2\lambda_A(t-u)}~I
  $$
In addition, when $\nabla_x b_t$ is uniformly bounded,  $\nabla_x \sigma_{k,t}=0$ and $  (H_C)$ is met, using (\ref{diff-contraction-general-pre}) we have almost sure estimate
   \begin{eqnarray*}
\Vert \partial_{\epsilon}Y^{\epsilon}_{s,t}\Vert&\leq& e^{-\lambda_A(t-s)} \Vert \overline{Z}_1-\overline{Z}_0\Vert+\Vert \nabla_x b_t\Vert_2~\int_s^t~e^{-\lambda_A(t-u)}~\EE\left(\Vert \partial_{\epsilon}X^{\epsilon}_{s,u}\Vert\right)~du\\
&\leq &e^{-\lambda_A(t-s)} \Vert \overline{Z}_1-\overline{Z}_0\Vert+\frac{\Vert \nabla_x b_t\Vert_2}{\lambda_A-\lambda_C}\left(e^{-\lambda_C(t-s)}-e^{-\lambda_A(t-s)}\right)~\EE\left(\Vert \overline{Z}_1-\overline{Z}_0\Vert^2\right)^{1/2}
    \end{eqnarray*}
    with the uniform spectral norm 
    $$
     \Vert \nabla_x b_t\Vert_2:=\sup_{x,y}\Vert \nabla_x b_t(x,y)\Vert_2
    $$
We summarize the above discussion with the following theorem.
\begin{theo}\label{theo-as-ref}
Assume that $\nabla_x b_t$ is uniformly bounded, $\nabla_x \sigma_{k,t}=0=\nabla_y \sigma_{k,t}$ and conditions $  (H_A)$ and $  (H_C)$ are met.
In this situation, we have the almost sure estimate
$$
\Vert \partial_{\epsilon}X^{\epsilon}_{s,t}\Vert\leq e^{-\lambda_A(t-s)} \Vert Z_1-Z_0\Vert+ (t-s)~e^{-\lambda(t-s)}~\Vert \nabla_x b_t\Vert_2~\EE\left(\Vert Z_1-Z_0\Vert^2\right)^{1/2}
$$
 with the process  $X^{\epsilon}_{s,t}$ defined in  (\ref{def-mu-epsilon}) and the parameter $\lambda:=\lambda_A\wedge \lambda_C$.
\end{theo}

\section{Some extensions}
    
 \subsection{A backward variational formula}

The stochastic transition semigroup associated with the flow $X^{\mu}_{s,t}(x)$ is defined for any mesurable function $f$
on $\RR^d$ by the formula
$$
\PP^{\mu}_{s,t}(f)(x):=f(X^{\mu}_{s,t}(x))\Longrightarrow P^{\mu}_{s,t}(f)(x)=\EE\left(\PP^{\mu}_{s,t}(f)(x)\right)
$$
For twice differentiable function $f$  we have the gradient and the Hessian formulae
\begin{eqnarray*}
\nabla \,\PP^{\mu}_{s,t}(f)(x)&=&\nabla X^{\mu}_{s,t}(x)~\PP^{\mu}_{s,t}(\nabla f)(x)\\
\nabla^2 \PP^{\mu}_{s,t}(f)(x)
&=&\left[\nabla X^{\mu}_{s,t}(x)\otimes \nabla X^{\mu}_{s,t}(x)\right]~\PP^{\mu}_{s,t}(\nabla^2 f)(x)+\nabla^2 X^{\mu}_{s,t}(x)~\PP^{\mu}_{s,t}(\nabla f)(x)
\end{eqnarray*}
In the above display,  $\nabla^2 X^{\mu}_{s,t}(x)$ stand for the tensors functions
\begin{eqnarray*}
\nabla^2\, X^{\mu}_{s,t}(x)_{(i,j),k}&=&\partial_{i,j}X^{\mu,k}_{s,t}(x)\\
\left[\nabla X^{\mu}_{s,t}(x)\otimes \nabla X^{\mu}_{s,t}(x)\right]_{(i,j),(k,l)}&=&\nabla X^{\mu}_{s,t}(x)_{i,k}\nabla X^{\mu}_{s,t}(x)_{j,l}
\end{eqnarray*}

Also recall that the infinitesimal generator $ L_{t,\phi_{s,t}(\mu)}$ of the stochastic flow (\ref{diff-st-ref-general}) is given for any twice differentiable function $f$ by the second order operator
      $$
 L_{t,\phi_{s,t}(\mu)}(f)(x):=  \left\langle b_{t}(\phi_{s,t}(\mu),x), \nabla f(x)\right\rangle+\frac{1}{2}~\tr\left[\nabla^2 f(x) ~\sigma_{t}(\phi_{s,t}(\mu),x)~\sigma_{t}(\phi_{s,t}(\mu),x)^{\prime}\right]
   $$
   Next theorem is an extension of a theorem by Da Prato-Menaldi-Tubaro~\cite{daprato-2} to nonlinear diffusions.
   \begin{theo}\label{back-theo-ref}
      Assume that $b_t(x,y)$ and $ \sigma_t(x,y)$ are Lipschitz functions w.r.t. the parameters $(t,x,y)$. 
   In this situation,  for any $\mu\in \Pa_2(\RR^d)$ we have
\begin{eqnarray}\label{backward-mck-v}
\PP^{\mu}_{s,t}(f)(x)
&=&f(x)+\int_s^t~L_{u,\phi_{s,u}(\mu)} \left(\PP^{\phi_{s,u}(\mu)}_{u,t}(f)\right)(x)~du\nonumber\\
&&\hskip4cm+\int_s^t~\nabla\, \PP^{\phi_{s,u}(\mu)}_{u,t}(f)(x)^{\prime}\,\sigma_{u}(\phi_{s,u}(\mu),x)\, \widehat{d}\,W_u
\end{eqnarray}
where $\widehat{d}\,W_u$ stands for the backward integration notation, so that the r.h.s. term in the above formula is a square integrable backward martingale. 
\end{theo}
The proof of the above formula follows the elegant stochastic backward 
variational analysis developed in~\cite{daprato-2}. A sketched proof is provided in the appendix, on page~\pageref{backward-mck-v-proof}.
 
We further assume that $ \nabla_x \sigma_{k,t}(x,y)=0$. In this situation, using the backward formula (\ref{backward-mck-v}) we check the stochastic interpolation formula
$$
\partial_u\left(X^{\phi_{s,u}(\mu)}_{u,t}\circ X^{\eta}_{s,u}\right)(y)^{\prime}=\left[\phi_{s,u}(\eta)-\phi_{s,u}(\mu)\right](b_u(\point,X^{\eta}_{s,u}(y)))^{\prime}\left[\nabla X^{\phi_{s,u}(\mu)}_{u,t}\right]({X}^{\eta}_{s,u}(y))
$$
Equivalently, we have
\begin{equation}\label{stoch-interpolation-bis}
 X^{\eta}_{s,t}(x)-X^{\mu}_{s,t}(x)=\int_s^t\left[\nabla X^{\phi_{s,u}(\mu)}_{u,t}\right]({X}^{\eta}_{s,u}(x))^{\prime}~\left[\phi_{s,u}(\eta)-\phi_{s,u}(\mu)\right](b_u(\point,X^{\eta}_{s,u}(x)))
~du
\end{equation}

Combining (\ref{ref-nablax-estimate-0}) and  (\ref{ref-nablax-estimate-0-ae}) with (\ref{diff-contraction-general}) we obtain the following corollary.

\begin{cor}\label{ref-intro-theo}
Assume the conditions of theorem~\ref{back-theo-ref} are satisfied and we have $ \nabla_x \sigma_{k,t}=0$ and $\Vert \nabla_x b_t(x,y)\Vert_2\leq c$, for some constant $c<\infty$.
Also assume that $(H_A)$ and $(H_C)$ are met  for some parameters $\lambda_A$ and  $\lambda_C$.  In this situation we have
 the exponential
decay estimates
$$
\EE\left(\Vert X^{\eta}_{s,t}(x)-X^{\mu}_{s,t}(x)\Vert^2\right)^{1/2}\leq c~\sqrt{d}~(t-s)~e^{-\lambda (t-s)}~\WW_2(\eta,\mu)\quad \mbox{\rm with}\quad \lambda:=\lambda_A\wedge \lambda_C
$$
In addition, when $\nabla_y \sigma_{k,t}=0$ we have the uniform and 
almost sure  estimates
$$
\Vert X^{\eta}_{s,t}(x)-X^{\mu}_{s,t}(x)\Vert\leq c~(t-s)~e^{-\lambda(t-s)}~\WW_2(\eta,\mu)
$$
\end{cor}

 \subsection{Diffusions on smooth manifolds}\label{sec-manifold}
 
 This section is concerned with the extension of our results to nonlinear diffusions on Riemannian manifolds. 
 Let us begin with the general necessary facts about nonlinear diffusions in manifolds. Our presentation will be made as similar as possible to the one in Euclidean space. For this, we will need It\^o differentials of manifold valued diffusions, parallel translation, covariant differential of tangent bundle valued semimartingales. 
 
 Let $M$ be a smooth manifold of dimension $d$. Stratonovich calculus is similar on $M$ and on $\RR^d$. So we are able to deal with Stratonovich SDE's of the type 
 \begin{equation}
  \label{Strato-SDE}
  \circ d X_{s,t}^\mu(x)= b^S_t(\phi_{s,t}(\mu), X_{s,t}^\mu(x))~dt+\sigma(X_{s,t}^\mu(x)))~\circ dW_t,
 \end{equation}
 where  for $y\in M$
 $$
 b^S_t(\eta,y)=\int_M\eta(dx)b^S_t(x,y),\quad b^S_t(x,y)\in T_yM,
 $$
 $W_t$ is a $\RR^m$-valued Brownian motion and $\sigma(y)$ is a linear map $\RR^m\to T_yM$. For simplicity $\sigma$ will not depend on time, but the time-dependent $\sigma$ can also be treated, we refer to~\cite{ACT10} for this extension, and also for the details of the constructions below. 

The only situation we will be interested in is when for all $y\in M$ the map 
$$
(\sigma\sigma^\ast)(y) : T_y^\ast M\to T_yM
$$
is a linear diffeomorphism. In this situation a scalar product can be defined in $T_y^\ast M$ and then in $T_yM$, leading to a Riemannian structure on $M$. The scalar product in $T_y^\ast M$ is 
\begin{equation}
 \label{spTastM}
 g^\ast(y)(\alpha,\beta)=\langle \sigma^\ast(y)(\alpha),\sigma^\ast(y)(\beta)\rangle_{\RR^m},
\end{equation}
and the scalar product in $T_yM$ is 
\begin{equation}
 \label{gyuv}
 g(y)(u,v)=g^\ast(y)\left((\sigma\sigma^\ast)^{-1}(y)(u),(\sigma\sigma^\ast)^{-1}(y)(v)\right).
\end{equation}
Associated to the metric $g$ is the Levi-Civita connection $\nabla$, which will be used to define parallel transport, It\^o equations, It\^o covariant differentials.
Recall that the parallel transport along a continuous $M$-valued semimartingale $X$ is the linear map 
$
//_t : T_{X_0}M\to T_{X_t}M 
$
which satisfies $//_0={\rm Id}$ and the Stratonovich SDE $\nabla_{\circ dX_t}//_t=0$. It is the natural extension to parallel transport along smooth paths, and it is an isometry. Parallel translation allows to anti-develop $X_t$ in $T_{X_0}M$ with the Stratonovich integral
$$
\As(X)_t=\int_0^t //_s^{-1}\circ dX_s
$$
The process $\Aa(X)$ takes its values in the vector space, it has an It\^o differential $d\Aa(X)_t$, which allows to define the It\^o differential of $X_t$
\begin{equation}
 \label{ItoXt}
 d^\nabla X_t:=//_td\As(X)_t.
\end{equation}
This It\^o differential is formally a vector which can be expressed in local coordinates as 
$$
d^\nabla X_t=\left(dX_t^i+\frac12~\Gamma_{j,k}^i(X_t)~d{<}X^j,X^k{>}_t\right) \frac{\partial}{\partial x^i}(X_t), \quad \hbox{with the Christoffel symbols }\quad \Gamma_{j,k}^i.
$$
The next object to consider is It\^o covariant derivative $DU_t$ of a $T_{X_t}M$-valued continuous semimartingale $U_t$:
\begin{equation}
 \label{Icd}
 DU_t:=//_td\left(//_t^{-1}U_t\right),
\end{equation}
easily defined from the fact that $//_t^{-1}U_t$ is vector valued. From the isometry property of parallel translation we easily get the formula for $V_t$ another $T_{X_t}M$-valued semimartingale and $\langle\cdot ,\cdot\rangle :=g$, 
\begin{equation}
 \label{DUV}
 d\langle U_t,V_t\rangle=\langle DU_t, V_t\rangle+\langle U_t,DV_t\rangle +\langle DU_t, DV_t\rangle.
\end{equation}

Defining $\displaystyle b_t(x,y):=b^S_t(x,y)+\frac12\sum_{k=1}^m\nabla\sigma_k(\sigma_k(y))$ (where for two vector fields $A,B$, $\nabla A(B(y))$ denotes the covariant derivative of $A$ in the direction $B(y)$), it is well known that the Stratonovich SDE~\eqref{Strato-SDE} is equivalent to the It\^o SDEs
\begin{equation}
 \label{Ito-SDE}
 d^\nabla X_{s,t}^\mu(x)= b_t(\phi_{s,t}(\mu), X_{s,t}^\mu(x))~dt+\sigma(X_{s,t}^\mu(x))~ dW_t.
\end{equation}
A remarkable fact concerning this equation, is that whenever it exists, a solution to equation~\eqref{Ito-SDE} is a diffusion with nonlinear generator $L_{t,\phi_{s,t}}(\mu)$, where 
\begin{equation}
 \label{genX}
 L_{t,\eta}=\frac12~\Delta +\int_M\eta(dx)~ b_t(x,y).
\end{equation}
So we can consider that the starting point of our study is SDE~\eqref{Ito-SDE} in a Riemannian manifold $(M,g)$. 

Let us adapt the regularity conditions $(H_A)$ and $(H_C)$: 

Define $A^g_t(x,y):=\nabla_y b_t(x,y)+\nabla_y b_t(x,y)^{\prime}$, where $\nabla_y b_t(x,y)$ is the covariant derivative with respect to the variable $y$, it is a linear map from $T_yM$ into itself, and $\nabla_y b_t(x,y)^{\prime}$ is its adjoint with respect to the Riemannian metric. 

{\em $(H^g_{A})$ : There exists some $\lambda_A^g\in \RR$ such that for any $x,y\in M$ and $t\geq 0$ we have 
   \begin{equation}\label{ref-mat-Upsilon-g}
\qquad A^g_t(x,y)-{\rm Ric}(y)\leq -2\lambda_A^g~g(y)
 \end{equation} 
 }
 where $\rm Ric$ is the Ricci curvature tensor of $M$.

Let $B_t^g$  be as in~\eqref{def-B-Sigma} with gradient replaced by covariant derivative. 

Define $C_t^g(x,y):= \frac{1}{2}\left[B^g_{t}\left(x,y\right)+ B^g_{t}\left(x,y\right)^{\prime}\right]$.

{\em $(H^g_{C})$ : There exists some $\lambda_C^g\in \RR$ such that for any $x,y\in M$ and $t\geq 0$ we have 
 \begin{equation}\label{hyp-b-sigma-general-g}
C_t^g(x,y)-\frac12{\rm Ric}_{M\times M}(x,y)
\leq -\lambda_C^g~g_{M\times M}(x,y)
 \end{equation}}
where $g_{M\times M}(x,y)$,  ${\rm Ric}_{M\times M}(x,y)$ are the product metric and Ricci curvature on $M\times M$.

\begin{theo}\label{stab-nl-theo-g}
 We have the exponential expansion or contraction inequalities

 \begin{equation}\label{ref-eta-mu-cv-g}
  (H^g_{A})\quad\Longrightarrow\quad \WW_2\left(\eta_0 P^{\mu}_{s,t},\eta_1 P^{\mu}_{s,t}\right)\leq c~e^{-\lambda_A^g(t-s)}~  \WW_2\left(\eta_0 ,\eta_1 \right)
 \end{equation}
for some finite constant $c$. In addition, we have
           \begin{equation}\label{diff-contraction-general-g}
(H^g_{C})\quad\Longrightarrow\quad \WW_2\left(\phi_{s,t}(\mu_0),\phi_{s,t}(\mu_1)\right)\leq e^{-\lambda_C^g (t-s)}~\WW_2(\mu_0,\mu_1)
 \end{equation}

   \end{theo}
\begin{remark}
The results of Theorem~\ref{stab-nl-theo-g} still hold
when $\sigma=\sigma_t$ and  $g=g_t$ depend on time, one just has to replace in $(H^g_{A})$  $\rm Ric$ by ${\rm Ric}-\dot g$ and in $(H^g_{C})$  $\rm Ric_{M\times M}$ by ${\rm Ric_{M\times M}}-\dot g_{M\times M}$. 
\end{remark}
\proof
The proof of the first estimate is  similar to the proof of Theorem~4.1 in \cite{ACT10} (where time dependent metrics are considered), so we will omit it. The proof of the second one is a combination of this proof and to the one of Theorem~\ref{stab-nl-theo} in the present article. Let us go into the details. 

Let $Z_0$, $Z_1$ two random variables with values in $M$, and such that $(Z_0,Z_1)$ minimizes $\EE[d^2(Z_0,Z_1)]$ under the condition that  $Z_0$ has law $\mu_0$ and $Z_1$ has law $\mu_1$. For all $\omega$, let $\epsilon \mapsto Z_\epsilon (\omega)$ be a geodesic between $Z_0(\omega)$ and $Z_1(\omega)$. 

As in the proof of Theorem~\ref{stab-nl-theo}, let $Y_{s,s}^{\mu_0}(x)=x$ and  $t\in [s,\infty[\mapsto Y_{s,t}^{\mu_0}(x)$ solve the equation
$$
dY_{s,t}^{\mu_0}(x)=b_t(\phi_{s,t}(\mu_0), Y_{s,t}^{\mu_0}(x))~dt +\sigma(Y_{s,t}^{\mu_0}(x))~d\bar W_t
$$
where $\bar W_t$ is a $\RR^m$ valued Brownian motion independent of $W_t$. Let $(\bar Z_\epsilon)_{\epsilon\in [0,1]}$ be independent of $( Z_\epsilon)_{\epsilon\in [0,1]}$ with the same law, $Y_{s,s}^\epsilon = \bar Z_\epsilon $ and $Y_{s,t}^\epsilon$ the solution to the It\^o SDE
\begin{equation}
\label{sde-curve}
dY_{s,t}^\epsilon =\EE_X\left[b_t(X_{s,t}^\epsilon, Y_{s,t}^\epsilon)\right]\,dt + //_{s,t}^{\,0,\epsilon}\,\left(\sigma (Y_{s,t}^0)~d\bar W_t\right),
\end{equation}
where $\epsilon\mapsto //_{s,t}^{\,0,\epsilon}(\omega)$ is the parallel transport along the  $\epsilon \mapsto Y_{s,t}^\epsilon(\omega)$. Notice that $Y_{s,t}^0\equiv Y_{s,t}^{\mu_0}(\bar Z_0)$. 

The equation~\eqref{sde-curve} is not an SDE on the manifold $M$, it is an SDE on $C^1$  $M$-valued paths. Existence of solutions have been established in~\cite{ACT10}. The processes $t\mapsto Y_{s,t}^\epsilon$ are obtained one from the others by infinitesimal synchronious coupling, and it is the only construction where a.s. the paths $\epsilon\mapsto Y_{s,t}^\epsilon(\omega)$ has finite variation. Moreover, the derivatives of theses paths satisfy
\begin{equation}
\label{sde-curve-derivative}
D\partial_\epsilon Y_{s,t}^\epsilon =\EE_X\left[\nabla_xb_t(X_{s,t}^\epsilon, Y_{s,t}^\epsilon)\partial_\epsilon X_{s,t}^\epsilon\right]\,dt + \EE_X\left[\nabla_yb_t(X_{s,t}^\epsilon, Y_{s,t}^\epsilon)\right]\partial_\epsilon Y_{s,t}^\epsilon\,dt-\frac12~ {\rm Ric}^\sharp (\partial_\epsilon Y_{s,t}^\epsilon)\,dt
\end{equation}
where ${\rm Ric}^\sharp(u)$ is the vector such that $\langle {\rm Ric}^\sharp(u),v\rangle ={\rm Ric}(u,v)$. The advantage of this construction is that the above covariant derivative has finite variation, and this implies 
$$
d\|\partial_\epsilon Y_{s,t}^\epsilon\|^2=2~\langle \partial_\epsilon Y_{s,t}^\epsilon, D\partial_\epsilon Y_{s,t}^\epsilon\rangle. 
$$ 
Then the proof is  similar to the one of Theorem~\ref{stab-nl-theo}: 
\begin{align*}
&\partial_t\,\EE\left[\|\partial_\epsilon Y_{s,t}^\epsilon\|^2\right]\\&=\EE\left[\left\langle \left(\begin{array}{c}\partial_\epsilon X_{s,t}^\epsilon\\\partial_\epsilon Y_{s,t}^\epsilon\end{array}\right), B_t(X^{\epsilon}_{s,t},Y^{\epsilon}_{s,t}) \left(\begin{array}{c}\partial_\epsilon X_{s,t}^\epsilon\\\partial_\epsilon Y_{s,t}^\epsilon\end{array}\right) \right\rangle\right] -\EE\left[{\rm Ric}(\partial_\epsilon Y_{s,t}^\epsilon,\partial_\epsilon Y_{s,t}^\epsilon)\right]\\
&= \EE\left[\left\langle \left(\begin{array}{c}\partial_\epsilon X_{s,t}^\epsilon\\\partial_\epsilon Y_{s,t}^\epsilon\end{array}\right), B_t(X^{\epsilon}_{s,t},Y^{\epsilon}_{s,t}) \left(\begin{array}{c}\partial_\epsilon X_{s,t}^\epsilon\\\partial_\epsilon Y_{s,t}^\epsilon\end{array}\right) \right\rangle\right] -\frac12\,\EE\left[{\rm Ric}_{M\times M}\left(\left(\begin{array}{c}\partial_\epsilon X_{s,t}^\epsilon\\\partial_\epsilon Y_{s,t}^\epsilon\end{array}\right),\left(\begin{array}{c}\partial_\epsilon X_{s,t}^\epsilon\\\partial_\epsilon Y_{s,t}^\epsilon\end{array}\right)\right)\right]\\
&\le -\lambda_C^g~\EE\left[ \left\|  \left(\begin{array}{c}\partial_\epsilon X_{s,t}^\epsilon\\\partial_\epsilon Y_{s,t}^\epsilon\end{array}\right)  \right\|^2\right]=-2\lambda_C^g\,\EE\left[\|\partial_\epsilon Y_{s,t}^\epsilon\|^2\right].
\end{align*}
 This implies that 
\begin{align*}
\EE\left[\|\partial_\epsilon Y_{s,t}^\epsilon\|^2\right]&\le \EE\left[\left\|\partial_\epsilon|_{\epsilon=0} Z_\epsilon\right\|^2\right]e^{-2\lambda_C^g(t-s)}=e^{-2\lambda_C^g(t-s)}~\WW_2^2(\mu_0,\mu_1).
\end{align*}
On the other hand, we have 
\begin{align*}
\WW_2^2\left(\phi_{s,t}(\mu_0),\phi_{s,t}(\mu_1)\right)&\le \EE\left[\left(\int_0^1\|\partial_\epsilon Y_{s,t}^\epsilon\|\, d\epsilon \right)^2\right]\\
&\le \int_0^1 \EE\left[ \|\partial_\epsilon Y_{s,t}^\epsilon\|^2 \right]\,d\epsilon \le ~e^{-2\lambda_C^g(t-s)}~\WW_2^2(\mu_0,\mu_1)
\end{align*}
This ends the proof of the theorem.
\cqfd

An important example of nonlinear diffusions in manifolds is again given by Langevin diffusions, defined as in~\eqref{Ito-SDE}, with now 
\begin{equation}
\label{E3.17}
b_t(x,y)=-\nabla U(y) - \nabla (F\circ\rho_x)(y)
\end{equation}
where $U$ is a potential function, $\rho$ is the Riemannian distance associated to the metric~$g$, $\rho_x$  is the distance to $x$ and $F : \RR_+\to \RR$ is a $C^2$ function. A sufficient condition $b_t(x,y)$ defined by~\eqref{E3.17} to be well defined and smooth is that the derivative of $F$ is nul at the origin and the support of $F$ is included in $[0,\imath(M))$, where   $\imath(M)$ denotes the injectivity radius of $M$. But smoothness of $b_t(x,y)$ is not a necessary condition for defining nonlinear diffusions. 

 We find that for $u,v\in T_yM$, 
\begin{equation}
\label{E3.18}
\nabla_y b (u,v) =-\nabla^2U (u,v)-\nabla^2(F\circ\rho_x)(u,v).
\end{equation}
In this context, condition $(H^g_{A}) $ reduces to 
\begin{equation}
\label{E3.19}
 \nabla^2U(y) +\nabla^2(F\circ\rho_x)(y)+\frac12~{\rm Ric}(y)\ge \lambda_A^g~ g(y).
\end{equation}
If for instance  $M$ is simply connected with nonpositive curvature (which implies that the distance function~$\rho$ is convex), and $F$ is nondecreasing, a sufficient condition is 
\begin{equation}
\label{E3.19ter}
 \nabla^2U(y) +\frac12~{\rm Ric}(y)\ge \lambda_A^g~ g(y).
\end{equation}

The computation of $B_t$ reveals that it is symmetric, and that for $(u,v)\in T_xM\times T_yM$, 
\begin{equation}
\label{E3.20}
\begin{split}
B_t(x,y)((u,v),(u,v))=&-\nabla^2U(x)(u,u)-\nabla^2U(y)(v,v)-\nabla^2(F\circ\rho)(x,y)((u,v),(u,v)),
\end{split}
\end{equation}
In this context condition $(H^g_{C}) $ reduces to
\begin{equation}
\label{E3.21}
\nabla^2U^{\oplus 2}(x,y)+\nabla^2(F\circ\rho)(x,y)+\frac12~{\rm Ric}_{M\times M}(x,y) \ge \lambda_C^g~ g_{M\times M}(x,y)
\end{equation}
where $U^{\oplus 2}(x,y)=U(x)+U(y)$.
Here again, when $M$ is simply connected with nonpositive curvature, $F$ is convex and nondecreasing, the above condition is met as soon as
\begin{equation}
\label{E3.21ter}
 \nabla^2U(y) +\frac12~{\rm Ric}(y)\ge \lambda_C^g ~ g(y).
\end{equation}
 \section{Mean field interacting diffusions}
 
 \subsection{Stability properties}\label{stab-particle-sec}
The interacting diffusion flow $\xi_{s,t}^j(z)=(\xi_{s,t}^{j,k}(z))_{1\leq k\leq d}\in\RR^d$ presented in (\ref{diff-st-ref-general-mf}) can be rewritten as
$$
 d\xi_{s,t}^{j}(z)=\Fa_t^{j}\left(\xi_{s,t}(z)\right)~dt+\sum_{1\leq \alpha\leq r} \Ga_{t,\alpha}^{j}\left(\xi_{s,t}(z)\right)~ dW^{j,\alpha}_t 
$$
with the drift and the diffusion functions defined for any $z=(z_1,\ldots,z_N)\in (\RR^d)^N$ with $z_i=(z_i^l)_{1\leq l\leq d}\in \RR^d$ by the formulae
$$
\Fa_t^{j,k}(z)=\frac{1}{N}~\sum_{1\leq n\leq N}~b^k_t(z_n,z_j)\quad \mbox{\rm and}\quad
\Ga_{t,\alpha}^{j,k}(z)=\frac{1}{N}~\sum_{1\leq n\leq N}~\sigma^k_{t,\alpha}(z_n,z_j)
$$
For any differentiable function $\Ha:z\in(\RR^d)^N\mapsto \Ha(z)\in (\RR^d)^N$ and any 
$1\leq i,j\leq N$ and $1\leq l,k\leq d$ we  consider the gradient blocks
$$
\left[\nabla\Ha (z)\right]_{i,j}=\nabla_{z_i}\Ha^j(z)\quad \mbox{\rm with}\quad
\left[\nabla_{z_i}\Ha^{j}(z)\right]_{l,k}=\partial_{z^l_i}\Ha^{j,k}(z)
$$
In this notation, for any $i\not=j$ we have
$$
\left[\nabla_{z_i}\Fa_t^{j}(z)\right]_{l,k}=
\frac{1}{N}~\partial_{z^l_i}b_t^k(z_i,z_j)\Longrightarrow \left[\nabla\Fa_t (z)\right]_{i,j}=\frac{1}{N}~\nabla_xb_t(z_i,z_j)
$$
and the diagonal term
$$
 \left[\nabla\Fa_t(z)\right]_{i,i}=\frac{1}{N}~\nabla_xb_t(z_i,z_i)+
\nabla_yb_t(m(z),z_i)
$$

Using the composition rule
\begin{equation}\label{c-rule}
 \nabla\left[\Ha_1\circ\Ha_2\right](z)=\nabla\Ha_2(z)~\left(\nabla\Ha_1\right)(\Ha_2(z))
\end{equation}
we check that
\begin{equation}\label{c-rule-cs}
 d\left[\nabla\xi_{s,t}(z)\right]_{i,j}=\left[\nabla\xi_{s,t}(z)\nabla\Fa_t\left(\xi_{s,t}(z)\right)\right]_{i,j}~dt+\sum_{1\leq \alpha\leq d} \left[\nabla\xi_{s,t}(z)\nabla\Ga_{t,\alpha}\left(\xi_{s,t}(z)\right)\right]_{i,j}~ dW^{j,\alpha}_t 
\end{equation}

 {\em $(\Ha_{\Aa})$ : There exists some $\lambda_{\Aa}\in \RR$ such that for any $z\in (\RR^d)^N$ and $t\geq 0$ we have 
   \begin{equation}\label{ref-mat-Upsilon-Aa}
\qquad \Aa_t(z)=\nabla\Fa_t(z)+\nabla\Fa_t(z)^{\prime}+\sum_{1\leq \alpha\leq r}\nabla\Ga_{t,\alpha}(z)
\nabla\Ga_{t,\alpha}(z)^{\prime}\leq -2\lambda_{\Aa}~I
 \end{equation} 
 }

 This spectral condition produces several gradient estimates. For instance, arguing as in (\ref{ref-nablax-estimate-0-ae-bis})  we have the following theorem.
 \begin{theo}\label{theo-Aa}
 Assume condition $(\Ha_{\Aa})$ is satisfied. In this situation we have the uniform exponential decay estimates
\begin{equation}\label{Ha-ref}
\EE\left[\Vert\nabla\xi_{s,t}(z)\Vert_2^2\right]^{1/2}\leq \EE\left[\Vert\nabla\xi_{s,t}(z)\Vert_F^2\right]^{1/2}\leq~\sqrt{dN}~e^{-\lambda_{\Aa} (t-s)} 
\end{equation}
 In addition, when $\nabla\Ga_{t,\alpha}(z)=0$ we have the uniform almost sure exponential decay estimate
 \begin{equation}
\Vert\nabla\xi_{s,t}(z)\Vert_2\leq e^{-\lambda_{\Aa} (t-s)} \quad \mbox{and}\quad
\Vert \xi_{s,t}(z)-\xi_{s,t}(\overline{z})\Vert\leq e^{-\lambda_{\Aa}(t-s)}~~\Vert z-\overline{z}\Vert\label{Ha-ref-2}
\end{equation}
\end{theo}
The proof of the above theorem is provided in the appendix, on page~\pageref{c-rule-proof}. 

  For the nonlinear Langevin diffusion discussed in (\ref{ex-langevin}) we have  $\nabla\Ga_{t,\alpha}(z)=0$ and
 \begin{eqnarray*}
 \left[\nabla\Fa_t (z)\right]_{i,j}&=&\frac{1}{N}~\nabla^2V(z_j-z_i)\\
  \left[\nabla\Fa_t(z)\right]_{i,i}&=&-\nabla^2 U(z_i)+\frac{1}{N}~\nabla^2V(0)-\frac{1}{N}\sum_{1\leq n\leq N}\nabla^2 V(z_i-z_n)
 \end{eqnarray*}
 In this situation we have
 $$
2^{-1} \Aa_t(z)=-\mbox{\rm Diag}\left(\nabla^2 U(z_1),\ldots,\nabla^2 U(z_N)\right)-\frac{1}{N}~\Ea_t(z)
 $$
 with the matrix $\Ea_t(z)$ with block entries
$$
 \left[\Ea_t(z)\right]_{i,j}=-\frac{1}{2}~\left[\nabla^2V(z_j-z_i)+\nabla^2V(z_i-z_j)\right]\quad \mbox{\rm and}\quad
  \left[\Ea_t(z)\right]_{i,i}=\sum_{1\leq n\not=i\leq N}\nabla^2 V(z_i-z_n)
$$
When $V$ is odd we have 
$$
\nabla^2 U(y)+\left(1-\frac{1}{N}\right)\nabla^2 V(y-x)\geq \lambda_{\Aa} I\Longrightarrow (\Ha_{\Aa})
$$
When $V$ is even and convex we have $\Ea_t(z)\geq 0$ and therefore
 $$
 \nabla^2 U(y)\geq \lambda_{\Aa}~I \Longrightarrow (\Ha_{\Aa})
 $$
In this situation, we also have
   \begin{equation}\label{ref-mean-field-Langevin}
 d\xi_t^i=-\frac{1}{N}\sum_{1\leq j\leq N}\left[\nabla V(\xi^i_t-\xi^j_t)+\nabla U(\xi^i_t)\right]~dt+ dW^i_t\quad \mbox{\rm with}\quad 1\leq i\leq N
 \end{equation}
Last but not least, whenever $\nabla V(0)=0$ we have 
 $$
 \begin{array}{l}
\displaystyle \Va(z):=\frac{1}{N}\sum_{1\leq i<j\leq N}
\frac{V(z_i-z_j)+V(z_j-z_i)}{2}+\sum_{1\leq i\leq N}U(z_i)\\
 \\
\displaystyle \Longrightarrow \nabla_{z_i}\Va (z)=\frac{1}{N}\sum_{1\leq j\leq N}\nabla V(z_i-z_j)+\nabla U(z_i)
 \end{array}
 $$
 Note that $\nabla V(0)=0$ holds when $V$ is even.
 In this situation, the diffusion $\xi_t$ reduces to a conventional Langevin diffusion 
 $$
  d\xi_t=-\nabla \Va(\xi_t)~dt+ d\Wa_t\quad \mbox{\rm with}\quad \Wa_t=(W_t^1,\ldots,W^N_t)^{\prime}
 $$
In this context, the stationary measure of the particle model $\xi_t$ is given by the Gibbs measure
 $$
 \nu(dz)\,\propto\, \exp{\left[-2\Va(z)\right]}~dz
 $$

 \subsection{Propagation of chaos properties}\label{prop-chaos-sec}
 
 For any differentiable function $g(x,y)$ from $\RR^{2d}$ into $\RR^d$ we let $\nabla_ug(x,y)$ be   the gradient matrices 
w.r.t.
 the coordinate $u\in\{x,y\}$, and we set
 $$
   \nabla_{x/y}g:= \nabla_{x}g+ \nabla_{y}g
 $$ 
We extend matrix-valued functions $G:z\in \RR^k\mapsto G(z)\in \RR^{d\times d}$ to the product space 
$\RR^{2k}$ by setting
$$
G[z;\overline{z}]:=\int_0^1G(\overline{z}+\epsilon (z-\overline{z}))~d\epsilon\quad\Longrightarrow\quad G[z;z]=G(z)
$$
We also consider the mapping $\delta : \RR^d\to  \RR^d\times  \RR^d$, $x\mapsto (x,x)$, and for any $x;\overline{x}\in\RR^d$ we set
$$
\nabla_{x/y}^\delta b_t[x;\overline{x}]\,:=\nabla_{x/y}(b_t\circ \delta)[x;\overline{x}]\quad \mbox{\rm and}\quad  
 \nabla_{x/y}^{\delta}\sigma_{t,k}[x;\overline{x}]\,:=\nabla_{x/y}(\sigma_{t,k}\circ \delta)[x;\overline{x}]$$
Let $  \Ba_t(z,\overline{z})$ and $  \Da_t(z,\overline{z})$ be the functions defined for any $z=(x,y)$ and $\overline{z}:=(\overline{x},\overline{y})\in\RR^{2d}$ by
       \begin{eqnarray*}
    \Ba_t(z,\overline{z})&:=&   \frac{1}{N} ~\Ba^{(1)}_t(z,\overline{z})+\left(1-\frac{1}{N}\right)~  \Ba^{(0)}_t(z,\overline{z})\\
     \Da_t(z,\overline{z})&:= & \frac{1}{N}~    \Da^{(1)}_t(z,\overline{z})+\left(1- \frac{1}{N}\right)~   \Da^{(0)}_t(z,\overline{z})
     \end{eqnarray*}  
The matrices $ \Ba^{(i)}_t(z,\overline{z})$ in the above display are given by
      \begin{eqnarray*}
    \Ba^{(1)}_t(z,\overline{z})&:=&\left[
\begin{array}{cc}
\nabla_{x/y}^{\delta} b_t[x;\overline{x}]&0\\
0&\nabla_{x/y}^{\delta}b_t[y;\overline{y}]
\end{array}
\right]~
\\
  \Ba^{(0)}_t(z,\overline{z})&:=&\left[
\begin{array}{cc}
\nabla_yb_t[(y,x);(\overline{y},\overline{x})]~&\nabla_xb_t[(y,x);(\overline{y},\overline{x})]\\
&\\
\nabla_xb_t[(x,y);(\overline{x},\overline{y}))]&\nabla_yb_t[(x,y);(\overline{x},\overline{y})]
\end{array}
\right]~
    \end{eqnarray*}
and the matrices $ \Da^{(i)}_t(z,\overline{z})$ are given by
         \begin{eqnarray*}
    \Da^{(1)}_t(z,\overline{z})&:=& \sum_{1\leq k\leq r}~\left[
\begin{array}{cc}
 \nabla_{x/y}^{\delta}\sigma_{t,k}[x;\overline{x}]^{\prime}~ \nabla_{x/y}^{\delta}\sigma_{t,k}[x;\overline{x}]
&0\\
0& \nabla_{x/y}^{\delta}\sigma_{t,k}[y;\overline{y}]^{\prime}~ \nabla_{x/y}^{\delta}\sigma_{t,k}[y;\overline{y}]
\end{array}
\right]\\
     \Da^{(0)}_t(z,\overline{z})&:=&2\sum_{1\leq k\leq r} 
     \left(  \begin{array}{cc} 
      \nabla_x\sigma_{t,k}[z;\overline{z}]^{\prime}~\nabla_x\sigma_{t,k}[z;\overline{z}]&
      \nabla_x\sigma_{t,k}[z;\overline{z}]^{\prime}~\nabla_y\sigma_{t,k}[z;\overline{z}]\\
       \nabla_y\sigma_{t,k}[z;\overline{z}]^{\prime}~\nabla_x\sigma_{t,k}[z;\overline{z}]    &
      \nabla_y\sigma_{t,k}[z;\overline{z}]^{\prime}~\nabla_y\sigma_{t,k}[z;\overline{z}]  
       \end{array}\right)
      \end{eqnarray*} 

Consider the following regularity condition:
  
     {\em $(\Ha_{\Ca})$ : There exists some $\lambda_{\Ca}\in \RR$ such that for any $z,\overline{z}\in \RR^{2d}$ and $t\geq 0$ we have 
 \begin{equation}\label{hyp-b-sigma-gen}
\Ca_t(z,\overline{z}):= \frac{1}{2}\left[\Ba_{t}(z,\overline{z})+ \Ba_{t}(z,\overline{z})^{\prime}\right]+
\Da_t(z,\overline{z})\leq -\lambda_{\Ca}~I
 \end{equation}}

Let $\zeta_0=(\zeta_0^i)_{1\leq i\leq N}$ be $N$ independent copies of a random variable with distribution $\mu$ on $\RR^d$.
Let $\xi_{t}:=\xi_{0,t}(\zeta_0)$  and consider the diffusion processes $\zeta_t=(\zeta_t^i)_{1\leq i\leq N}$ defined as $\xi_{t}$ 
by replacing the occupation measures $m(\xi_t)$ by the distributions $ \mu_t=\phi_{t}(\mu):=\phi_{0,t}(\mu)$; that is, for any $1\leq i\leq N$ we have
 $$
   d\zeta_t^i=b_t(\mu_t,\zeta^i_t)~dt+\sigma_{t}(\mu_t,\zeta^i_t)~ dW^{i}_t   
    $$

     \begin{theo}
     \label{T4.2}
  Assume condition   $(\Ha_{\Ca})$ is satisfied.
In this situation, for any $\epsilon>0$ and any distribution $\mu$ on $\RR^d$ we have
\begin{equation}
 \label{unif-prop-chaos-moments}
 \EE\left(\Vert \xi_t^1-\zeta_t^1\Vert^2\right)\leq \frac{1}{N}~\int_0^t~e^{-2(\lambda_{\Ca}-\epsilon)(t-s)}~   \left(
      2\alpha_s(\mu)  +\frac{ \beta_s(\mu)}{2\epsilon}\right)~ds
   \end{equation}
      with the parameters
    \begin{eqnarray*}
    \alpha_t(\mu)&:=&\sum_{1\leq k\leq r}~\int~\phi_{t}(\mu)(dx)~\Vert
\sigma_{t,k}(x,x)-\sigma_{t,k}(\phi_{t}(\mu),x)\Vert^2\\
    \beta_t(\mu)&:=&\frac1N\int~\phi_{t}(\mu)(dx)~\Vert
b_t(x,x)-b_t(\phi_{t}(\mu),x)\Vert^2\\
&&\hskip3cm+\left(1-\frac1N\right)\int~\phi_{t}(\mu)(dx)\phi_{t}(\mu)(dy)~\Vert
b_t(x,y)-b_t(\phi_{t}(\mu),y)\Vert^2
      \end{eqnarray*}

 \end{theo}
 
 \proof
 We set $S_t:=\EE\left(\Vert \xi_t^1-\zeta_t^1\Vert^2\right)$.
 Using the decomposition
    $$
    \begin{array}{l}
   d(\xi_t^1-\zeta_t^1)=\left[b_t(m(\xi_t),\xi^1_t)-b_t(\mu_t,\zeta^1_t)\right]~dt  +\left[\sigma_{t}(m(\xi_t),\xi^1_t)-\sigma_{t}(\mu_t,\zeta^1_t)\right]~ dW^{1}_t
    \end{array}
    $$
we check that
       $$
    \begin{array}{l}
  \displaystyle \partial_tS_t=2~\EE\left(
   \langle \xi_t^1-\zeta_t^1,b_t(m(\xi_t),\xi^1_t)-b_t(\mu_t,\zeta^1_t)\rangle\right)+   \Sigma_t+ \Gamma_t\\
\\
\hskip3cm  \displaystyle+   2\sum_{1\leq k\leq r}  \EE\left(\left\langle \sigma_{t,k}(m(\xi_t),\xi^1_t)-\sigma_{t,k}(m(\zeta_t),\zeta^1_t),\sigma_{t,k}(m(\zeta_t),\zeta^1_t)-\sigma_{t,k}(\mu_t,\zeta^1_t)\right\rangle\right)   \end{array}
    $$
      with $\Sigma_t$ and $  \Gamma_t$ defined by
    \begin{eqnarray*}
    \Sigma_t&:=&  \sum_{1\leq k\leq r}\EE\left(\Vert\sigma_{t,k}(m(\xi_t),\xi^1_t)-\sigma_{t,k}(m(\zeta_t),\zeta^1_t)\Vert^2\right)\\
    \Gamma_t&:=&   \sum_r  \EE\left(\Vert\sigma_{t,k}(m(\zeta_t),\zeta^1_t)-\sigma_{t,k}(\mu_t,\zeta^1_t)\Vert^2\right)
    \end{eqnarray*}
Applying Cauchy-Schwartz inequality  we find that
      $$
  \displaystyle 2^{-1} \partial_tS_t\leq I_t+\Sigma_t+J_t+ \Gamma_t
    $$
  with
      \begin{eqnarray*}
 I_t&:=&   \EE\left(
   \langle \xi_t^1-\zeta_t^1,b_t(m(\xi_t),\xi^1_t)-b_t(m(\zeta_t),\zeta^1_t)\rangle\right)\\
   J_t&:=&  \EE\left(\langle \xi_t^1-\zeta_t^1,b_t(m(\zeta_t),\zeta^1_t)-b_t(\mu_t,\zeta^1_t)\rangle\right)
    \end{eqnarray*}
  To estimate the term $\Sigma_t$
we observe that
    $$
       \begin{array}{l}
  \displaystyle
  \EE\left(\Vert\sigma_{t,k}(m(\xi_t),\xi^1_t)-\sigma_{t,k}(m(\zeta_t),\zeta^1_t)\Vert^2\right)\\
  \\
   \displaystyle =\frac{1}{N^2}\sum_{1\leq i,j\leq N} \EE\left(\left\langle\sigma_{t,k}(\xi_t^i,\xi^1_t)-\sigma_{t,k}(\zeta_t^i,\zeta^1_t),
\sigma_{t,k}(\xi_t^j,\xi^1_t)-\sigma_{t,k}(\zeta_t^j,\zeta^1_t)\right\rangle\right)\\
\\
   \displaystyle \leq \frac{1}{N^2}\sum_{1\leq i,j\leq N} \EE\left(\Vert\sigma_{t,k}(\xi_t^i,\xi^1_t)-\sigma_{t,k}(\zeta_t^i,\zeta^1_t)\Vert
\Vert\sigma_{t,k}(\xi_t^j,\xi^1_t)-\sigma_{t,k}(\zeta_t^j,\zeta^1_t)\Vert\right)\\
\\
   \displaystyle \leq \frac{1}{N}~\EE\left(\Vert\sigma_{t,k}(\xi_t^1,\xi^1_t)-\sigma_{t,k}(\zeta_t^1,\zeta^1_t)\Vert^2\right)  
   +\left(1- \frac{1}{N}\right)~\EE\left(\Vert\sigma_{t,k}(\xi_t^1,\xi^2_t)-\sigma_{t,k}(\zeta_t^1,\zeta^2_t)\Vert^2\right)    \end{array}     
    $$
    On the other hand, for any differentiable function $g$ from $\RR^{2d}$ into $\RR^d$, and for any $z=(x,y)$ and $\overline{z}=(\overline{x},\overline{y})\in \RR^{2d}$
 we have the first order decomposition
    $$
          \begin{array}{l}
  \displaystyle   
g(z)-g(\overline{z})= \nabla_xg[z,\overline{z}]~(x-\overline{x})+
  \nabla_yg[z,\overline{z}]~(y-\overline{y})\\
  \\
    \displaystyle   \Longrightarrow\quad   \left\{
           \begin{array}{rcl}
    \Vert g(z)-g(\overline{z})  \Vert^2
 & =& \left(  
x-\overline{x},y-\overline{y}\right)^{\prime}~ \overline{\nabla}g[z,\overline{z}]~\left(  \begin{array}{l}
      x-\overline{x}\\
      y-\overline{y}
          \end{array}\right) \\
&&\\
 \Vert g(x,x)-g(\overline{x},\overline{x})  \Vert^2&=&
 ~(x-\overline{x})^{\prime} ~  \nabla_{x/y}g[x,\overline{x}]^{\prime}~ \nabla_{x/y}g[x,\overline{x}]~(x-\overline{x})          
                      \end{array}\right.
    \end{array}
    $$ 
with the matrix
    $$
    \overline{\nabla}g[z,\overline{z}]:=\left(  \begin{array}{cc} 
      \nabla_xg[z,\overline{z}]^{\prime}~\nabla_xg[z,\overline{z}]&
      \nabla_xg[z,\overline{z}]^{\prime}~\nabla_yg[z,\overline{z}]\\
       \nabla_yg[z,\overline{z}]^{\prime}~\nabla_xg[z,\overline{z}]     &
      \nabla_yg[z,\overline{z}]^{\prime}~\nabla_yg[z,\overline{z}]  
       \end{array}\right)
    $$ 
By symmetry arguments, this implies that
     $$
                \begin{array}{l}
  \displaystyle 
2\sum_{1\leq k\leq r}\EE\left(\Vert\sigma_{t,k}(\xi_t^1,\xi^1_t)-\sigma_{t,k}(\zeta_t^1,\zeta^1_t)\Vert^2\right)\\
\\
 \displaystyle= \EE\left[\left(  
\xi_t^1-\zeta_t^1,\xi_t^2-\zeta_t^2\right)^{\prime}~\Da^{(1)}_t((\xi_t^1,\xi^2_t),(\zeta_t^1,\zeta^2_t))~\left(  \begin{array}{l}
      \xi_t^1-\zeta_t^1\\
   \xi_t^2-\zeta_t^2
          \end{array}\right)\right]
    \end{array}   $$ 
    In the same vein, we have
    $$
           \begin{array}{l}
  \displaystyle 2\sum_{1\leq k\leq r}
    ~\EE\left(\Vert\sigma_{t,k}(\xi_t^1,\xi^2_t)-\sigma_{t,k}(\zeta_t^1,\zeta^2_t)\Vert^2\right) \\
    \\
    \displaystyle  = \EE\left[\left(  
\xi_t^1-\zeta_t^1,\xi_t^2-\zeta_t^2\right)^{\prime}~   \Da^{(0)}_t((\xi_t^1,\xi^2_t),(\zeta_t^1,\zeta^2_t))~\left(  \begin{array}{l}
      \xi_t^1-\zeta_t^1\\
   \xi_t^2-\zeta_t^2
          \end{array}\right)\right]
  \end{array}     $$
   This yields the estimate
       \begin{equation}\label{ref-Sigmat}
  \displaystyle 2 \Sigma_t \leq ~ \EE\left[\left(  
\xi_t^1-\zeta_t^1,\xi_t^2-\zeta_t^2\right)^{\prime}~   \Da_t((\xi_t^1,\xi^2_t),(\zeta_t^1,\zeta^2_t))~\left(  \begin{array}{l}
      \xi_t^1-\zeta_t^1\\
   \xi_t^2-\zeta_t^2
          \end{array}\right)\right]
`\end{equation}
To estimate the term $ I_t$ we use the decomposition
       \begin{equation}
\label{ItRd}
I_t =\frac{1}{N}~\EE\left(
   \langle \xi_t^1-\zeta_t^1,b_t(\xi^1_t,\xi^1_t)-b_t(\zeta_t^1,\zeta^1_t)\rangle\right)+\left(1-\frac{1}{N}\right)~\EE\left(
   \langle \xi_t^1-\zeta_t^1,b_t(\xi^2_t,\xi^1_t)-b_t(\zeta_t^2,\zeta^1_t)\rangle\right)
    \end{equation}

    Also notice that
    $$
     \begin{array}{l}
  \displaystyle
2~ \EE\left(  \langle \xi_t^1-\zeta_t^1,b_t(\xi^1_t,\xi^1_t)-b_t(\zeta_t^1,\zeta^1_t)\rangle\right)\\
   \\
  \displaystyle   = 2~\EE\left( 
   \langle \xi_t^1-\zeta_t^1,\nabla_{x/y}(b_t\circ \delta)[\xi^1_t,\zeta_t^1]~(\xi^1_t-\zeta_t^1)
\rangle\right)\\
\\
=\EE\left[\left(  
\xi_t^1-\zeta_t^1,\xi_t^2-\zeta_t^2\right)^{\prime}~\Ba^{(1)}_t[(\xi^1_t,\xi^2_t);(\zeta_t^1,\zeta^2_t)]~\left(  \begin{array}{l}
      \xi_t^1-\zeta_t^1\\
   \xi_t^2-\zeta_t^2
          \end{array}\right)\right]
       \end{array}
    $$ 
We also have
       $$
     \begin{array}{l}
  \displaystyle
 \EE\left(  \langle \xi_t^1-\zeta_t^1,b_t(\xi^2_t,\xi^1_t)-b_t(\zeta_t^2,\zeta^1_t)\rangle\right)\\
   \\
  \displaystyle   =
   \EE\left( \langle \xi_t^1-\zeta_t^1,\nabla_xb_t[(\xi^2_t,\xi^1_t);(\zeta_t^2,\zeta^1_t)]~(\xi^2_t-\zeta_t^2)
\rangle\right)
+  \EE\left( \langle \xi_t^1-\zeta_t^1,
\nabla_yb_t[(\xi^2_t,\xi^1_t);(\zeta_t^2,\zeta^1_t)]~(\xi^1_t-\zeta_t^1)
\rangle\right)
       \end{array}
    $$ 
    This implies that
         $$
     \begin{array}{l}
  \displaystyle
2 \EE\left(  \langle \xi_t^1-\zeta_t^1,b_t(\xi^2_t,\xi^1_t)-b_t(\zeta_t^2,\zeta^1_t)\rangle\right)\\
   \\
  \displaystyle   =
   \EE\left( \langle \xi_t^1-\zeta_t^1,\nabla_xb_t[(\xi^2_t,\xi^1_t);(\zeta_t^2,\zeta^1_t)]~(\xi^2_t-\zeta_t^2)
\rangle\right)\\
\\
 \hskip1cm+  \EE\left( \langle \xi_t^2-\zeta_t^2,\nabla_xb_t[(\xi^1_t,\xi^2_t);(\zeta_t^1,\zeta^2_t)]~(\xi^1_t-\zeta_t^1)
\rangle\right)\\
\\
 \hskip2cm+  \EE\left( \langle \xi_t^1-\zeta_t^1,
\nabla_yb_t[(\xi^2_t,\xi^1_t);(\zeta_t^2,\zeta^1_t)]~(\xi^1_t-\zeta_t^1)
\rangle\right)\\
\\
 \hskip3cm+  \EE\left( \langle \xi_t^2-\zeta_t^2,
\nabla_yb_t[(\xi^1_t,\xi^2_t);(\zeta_t^1,\zeta^2_t)]~(\xi^2_t-\zeta_t^2)
\rangle\right)
       \end{array}
    $$ 
    from which we check that
      $$
     \begin{array}{l}
  \displaystyle
2\EE\left(  \langle \xi_t^1-\zeta_t^1,b_t(\xi^2_t,\xi^1_t)-b_t(\zeta_t^2,\zeta^1_t)\rangle\right)\\
   \\
  \displaystyle   =\EE\left[\left(  
\xi_t^1-\zeta_t^1,\xi_t^2-\zeta_t^2\right)^{\prime}~ \Ba^{(0)}_t[(\xi^1_t,\xi^2_t);(\zeta_t^1,\zeta^2_t)]~\left(  \begin{array}{l}
      \xi_t^1-\zeta_t^1\\
   \xi_t^2-\zeta_t^2
          \end{array}\right)\right]       \end{array}
    $$ 
 Combining the above decompositions we check that
      $$
 2 I_t =\EE\left[\left(  
\xi_t^1-\zeta_t^1,\xi_t^2-\zeta_t^2\right)^{\prime}~  \Ba_t[(\xi^1_t,\xi^2_t);(\zeta_t^1,\zeta^2_t)]~\left(  \begin{array}{l}
      \xi_t^1-\zeta_t^1\\
   \xi_t^2-\zeta_t^2
          \end{array}\right)\right]  
    $$  
Combining the above estimate with (\ref{ref-Sigmat}) we find that
     $$
 \partial_tS_t\leq \EE\left[\left(  
\xi_t^1-\zeta_t^1,\xi_t^2-\zeta_t^2\right)^{\prime}~  \Ca_t[(\xi^1_t,\xi^2_t);(\zeta_t^1,\zeta^2_t)]~\left(  \begin{array}{l}
      \xi_t^1-\zeta_t^1\\
   \xi_t^2-\zeta_t^2
          \end{array}\right)\right]   +2  J_t(\xi_t,\zeta_t)+2 \Gamma_t(\zeta_t)    
    $$
     from which we conclude that
      $$
     2^{-1} \partial_tS_t \leq-\lambda_{\Ca}~S_t   +  J_t+ \Gamma_t 
      $$
 Applying twice Cauchy-Schwartz inequality we check the estimate
      $$
 \vert  J_t(\xi_t,\zeta_t)\vert\leq  \sqrt{S_t}~\EE\left(\Vert b_t(m(\zeta_t),\zeta^1_t)-b_t(\mu_t,\zeta^1_t)\Vert^2\right)^{1/2}
    $$
    On the other hand, we have
       \begin{eqnarray*}
\EE\left(\Vert b_t(m(\zeta_t),\zeta^1_t)-b_t(\mu_t,\zeta^1_t)\Vert^2\right)
&=& \frac{1}{N}~    \beta_t(\mu)\quad \mbox{\rm and}\quad     \Gamma_t(\zeta_t) =\frac{1}{N}~\alpha_t(\mu)    \end{eqnarray*}
    This implies that
 $$
     2^{-1} \partial_tS_t\leq-\lambda_{\Ca}~S_t  +  \frac{1}{\sqrt{N}}~\sqrt{\beta_t(\mu)S_t}  +\frac{1}{N}~\alpha_t(\mu)    
      $$
      Recalling that $2ab\leq \epsilon a^2+b^2/\epsilon$ for any $\epsilon>0$  and $a,b\in \RR$, we check that
       $$
   \partial_tS_t\leq-2(\lambda_{\Ca}-\epsilon)~S_t  +\frac{1}{N}~\left(
   2\alpha_t(\mu)  +\frac{ \beta_t(\mu)}{2\epsilon}  \right)
      $$
   This ends the proof of the theorem. \cqfd
   
We end this section with some comments on the regularity condition  $(\Ha_{\Ca})$.

     For the nonlinear Langevin diffusion discussed in (\ref{ex-langevin}) we have $\Da_t(z,\overline{z})=0$ and
     \begin{eqnarray*}
  \nabla_{x}b_t[(x,y);(\overline{x},\overline{y})]&=& \nabla^2V[y-x;\overline{y}-\overline{x}]\\
  \nabla_{y}b_t[(x,y);(\overline{x},\overline{y})]&=&-\nabla^2U[y;\overline{y}]- \nabla^2V[y-x;\overline{y}-\overline{x}]\quad\mbox{\rm and}\quad
    \nabla_{x/y}^{\delta}b_t[y;\overline{y}]=-\nabla^2U[y;\overline{y}] \end{eqnarray*}
In this context, we have
$$
\begin{array}{l}
\displaystyle   - \Ca_t(z,\overline{z}):=\left[
\begin{array}{cc}
\nabla^2U[x;\overline{x}]&0\\
0&\nabla^2U[y;\overline{y}]
\end{array}
\right]\\
\\
\hskip2cm\displaystyle+\left(1-\frac{1}{N}\right)\left[
\begin{array}{cc}
 \nabla^2V[x-y;\overline{x}-\overline{y}]&-\frac{\nabla^2V[x-y;\overline{x}-\overline{y}]+\nabla^2V[y-x;\overline{y}-\overline{x}]}{2}\\
 &\\
-\frac{\nabla^2V[x-y;\overline{x}-\overline{y}]+\nabla^2V[y-x;\overline{y}-\overline{x}]}{2}& \nabla^2V[y-x;\overline{y}-\overline{x}]
\end{array}
\right]
\end{array}
$$
Also observe that for any $z\in \RR^{2d}$ we have the decomposition
$$
\Ca_t^{(1)}(z):=\Ca_t(z,z)= \left(1-\frac{1}{N}\right)~C_t(z)+\frac{1}{N}~ C^{(1)}_t(z)
$$
 with the matrices
 $$
 C^{(1)}_t(z):=\frac{1}{2}\left[B^{(1)}_{t}(z)+ B^{(1)}_{t}(z)^{\prime}\right]+
D^{(1)}_t(z)
 $$
In the above display $C_t(z)$ stands for the matrix defined in (\ref{hyp-b-sigma-general}),  $B^{(1)}_{t}(z)$ and $D^{(1)}_{t}(z)$ stand for the matrices defined for any $z=(x,y)\in \RR^{2d}$ by
  \begin{eqnarray*}
   B^{(1)}_t(x,y)&:=&\left[
\begin{array}{cc}
\nabla_{x/y} b_t(x,x)&0\\
0&\nabla_{x/y}b_t(y,y)
\end{array}
\right]~
\\
    D^{(1)}_t(x,y)&:=& \sum_{1\leq k\leq r}~\left[
\begin{array}{cc}
 \nabla_{x/y}\sigma_{t,k}(x,x)^{\prime}~ \nabla_{x/y}\sigma_{t,k}(x,x)
&0\\
0& \nabla_{x/y}\sigma_{t,k}(y,y)^{\prime}~ \nabla_{x/y}\sigma_{t,k}(y,y)
\end{array}
\right]
      \end{eqnarray*} 
      
Consider the following regularity condition:
  
     {\em $\left(\Ha_{\Ca^{(1)}}\right)$ : There exists some $\lambda_{\Ca^{(1)}}\in \RR$ such that for any $(x,y)\in \RR^{2d}$ and $t\geq 0$ we have 
 \begin{equation}\label{hyp-b-sigma-gen-1}
\Ca_t^{(1)}(x,y)\leq -\lambda_{\Ca^{(1)}}~I
 \end{equation}}
Assume that $\left(\Ha_{\Ca^{(1)}}\right)$ is met. Using  the fact that $ \EE(\Sigma^{\prime})\EE(\Sigma)\leq \EE(\Sigma^{\prime}\Sigma)$, for any random matrix $\Sigma$, we check that
 \begin{equation}\label{hyp-b-sigma-gen-simple}
(\Ha_{\Ca})\quad \mbox{\rm and}\quad
(\Ha_{C})\quad \mbox{\rm are met}\quad\quad \mbox{\rm with}\quad \lambda_{\Ca}= \lambda_{\Ca^{(1)}}\quad \mbox{\rm and}\quad \lambda_C= \left(1-\frac{1}{N}\right)^{-1}~\lambda_{\Ca^{(1)}}
 \end{equation}

  Several uniform estimates can be derived combining (\ref{unif-prop-chaos-moments}) with the moments estimates  (\ref{unif-moments}).
   For instance,  suppose we are given a time homogeneous model $(b_t,\sigma_t)=(b,\sigma)$, for some functions $(b,\sigma)$
   with uniformly bounded first order derivatives. Also assume 
   $\left(\Ha_{\Ca^{(1)}}\right)$ is met for some
 $\lambda_{\Ca^{(1)}}>0$. In this context, the moments estimates  (\ref{unif-moments}) ensure that
$$
\alpha_t(\mu)\vee  \beta_t(\mu)\leq c(\mu)
$$
for some constant $c(\mu)$ whose values only depends on the measure $\mu$. Choosing $\epsilon=\lambda_{\Ca}/2$ in (\ref{unif-prop-chaos-moments}) we readily check that
  $$
   \EE\left(\Vert \xi_t^1-\zeta_t^1\Vert^2\right)\leq \frac{c(\mu) }{N}~\frac{1}{\lambda_{\Ca}}~\left(
      2 +\frac{ 1}{\lambda_{\Ca}}\right)~
  $$ 
   
\subsection{Propagation of chaos in manifolds}
\label{subsec4.3}

Our aim is to state an analogous of Theorem~\ref{T4.2} in a Riemannian manifold $(M,g)$. We will take the notations of Section~\ref{sec-manifold}. Let us denote by $\rho$ the Riemannian distance in $M$. 
Now $\zeta_0=(\zeta_0^i)_{1\le i\le N}$ are independent copies of a random variable with distribution $\mu$ on $M$. For $1\le i\le N$ the diffusions $\zeta_{s,t}^i(x)$ satisfy the It\^o SDE 
\begin{equation}
 \label{E4.3.1}
 d^\nabla \zeta^i_{s,t}(x)=b_t(\phi_{s,t}(\mu), \zeta^i_{s,t}(x))\,dt +\sigma(\zeta_{s,t}^i(x))\, dW_t^i,
\end{equation}
with $\sigma(y) :\RR^m\to T_yM$ linear,  $\sigma \sigma^\ast=g^\ast$, and $(W_t^i), 1\le i \le N$ independent $\RR^m$-valued Brownian motions independent of $\zeta_0$. 
Denote $\mu_t:=\phi_{0,t}(\mu)$, $\zeta_t:=\zeta_{0,t}(\zeta_0)$. The diffusions $\zeta_t^i$ are independent and identically distributed, with law $\mu_t$ at time $t$. Define an approximation of $\zeta_t$ with the Markov process $\xi_t=(\xi_t^i)_{1\le i\le N}$ satisfying $\xi_0=\zeta_0$ and for all $i$,  
\begin{equation}
 \label{E4.3.2}
 d^\nabla \xi^i_t=b_t(m(\xi_t), \xi^i_t)\,dt +//_{\zeta_t^i,\xi_t^i}(\sigma(\zeta_{t}^i)\, dW_t^i)
\end{equation}
where for $x,y\in M$, $//_{x,y}$ denotes parallel translation along the minimal geodesic from $x$ to $y$. It is well-known that such an equation has a solution, which realizes the coupling by parallel translation of martingale parts of $\zeta_t^i$ and $\xi_t^i$ (see e.g. \cite{ATW06} or~\cite{WangFY14}). The only difficulty is when $\xi_t^i$ is in the cutlocus of $\zeta_t^i$, but this difficulty can be overcome by constructing approximations of the solutions which are decoupled in an $\epsilon$-neighbourhood of the cutlocus, and by letting then $\epsilon$ tend to $0$. However the solution obtained is not strong. Anyway, since $//_{\zeta_t^i,\xi_t^i}$ is an isometry and the $W_t^i$ are independent, the process $\xi_t$ is a Brownian motion in $M^N$ with drift $(b_t(m(\xi_t), \xi^i_t))_{1\le i\le N}$, so it is a diffusion process. Moreover independent $\RR^m$ valued Brownian motions $w_t^i$ can be found such that
\begin{equation}
 \label{E4.3.3}
 d^\nabla \xi^i_t=b_t(m(\xi_t), \xi^i_t)\,dt +\sigma(\xi_{t}^i)\, dw_t^i,
\end{equation}
they satisfy $$dw_t^i=\sigma^\ast(\xi_t^i)~ (\sigma\sigma^\ast)^{-1}(\xi_t^i)~ //_{\zeta_t^i,\xi_t^i}(\sigma(\zeta_{t}^i)\, dW_t^i)+ dm_t^i$$ for some ``complementary'' martingale $m_t^i$.

The important fact about this construction is that the distance $\rho^2(\zeta_t^i, \xi_t^i)$ has finite variation. More precisely, letting for $x,y\in M$ with $y$ not belonging to the cutlocus of $x$, $s\mapsto \gamma(x,y)(s)$ the geodesic from $x$ to $y$ in time $1$ and $\overrightarrow{xy}=\dot\gamma(x,y)(0)$  we have 
\begin{equation}
 \label{E4.3.4}
 \begin{split}
 d\rho^2(\zeta_t^1, \xi_t^1)&=2\langle\dot\gamma(\zeta_t^1, \xi_t^1)(1), b_t(m(\xi_t), \xi^i_t)\rangle \,dt
 -2\langle \dot\gamma(\zeta_t^1, \xi_t^1)(0), b_t(\mu_t,\zeta^i_t)\rangle\,dt \\&\hskip3cm+2\rho(\zeta_t^1, \xi_t^1)~\frac12I(\zeta_t^1, \xi_t^1)\,dt-dL_t
 \end{split}
\end{equation}
In the above display $L_t$ stands for a nondecreasing process which increases only when $\xi_t^1$ is in the cutlocus of $\zeta_t^1$, and $I$ is the index map defined for $x,y\in M$, and $y\not\in{\rm Cut(x)}$, by
\begin{equation}
 \label{E4.3.5}
I(x,y)=\sum_{i=1}^{d-1}\int_0^{\rho(x,y)}\left(\left\|\nabla_{\dot\varphi(s)}J_i(s)\right\|^2-\langle R(\dot\varphi(s), J_i(s))J_i(s),\dot\varphi(s)\rangle\right)\,ds 
\end{equation}
where $\varphi$ is a unit speed geodesic from $x$ to $y$ started at time $0$, $(J_i(0))_{1\le i\le d-1}$ is an orthonormal basis of $\dot\varphi(0)^\perp$, $J_i(\rho(x,y))=//_{x,y}J_i(0)$ and $s\mapsto J(s)$ is a Jacobi field along $s\mapsto\varphi(s)$ (see e.g. \cite{ATW06}). It is well known that when ${\rm Ric}_M\ge \kappa$ then $I(x,y)\le \bar I(\rho(x,y),\kappa)$ where $\bar I(\rho(x,y),\kappa)$ is the same quantity computed in a constant curvature manifold, for two points at the same distance. Moreover we have the explicit values 
\begin{equation}
 \label{E4.3.6}
 \bar I(\rho,\kappa)=\left\{
 \begin{array}{ccc}
 \displaystyle -2\sqrt{(d-1)\kappa}~\tan\left(\frac{\rho}{2}\sqrt{\frac{\kappa}{d-1}}\right)&\hbox{if}&\kappa > 0\\
  0&\hbox{if}&\kappa=0\\
   \displaystyle 2\sqrt{(d-1)(-\kappa)}~\tanh\left(\frac{\rho}{2}\sqrt{\frac{-\kappa}{d-1}}\right)&\hbox{if}&\kappa < 0
 \end{array}
 \right.
\end{equation}
In any case, $\bar I(\rho,\kappa)\le -\kappa\rho$, so we obtain as a general result that when ${\rm Ric}_M\ge \kappa$
\begin{equation}
 \label{E4.3.7}
 I(x,y)\le -\kappa\,\rho(x,y).
\end{equation}
So we have 
\begin{equation}
 \label{E4.3.8}
 \begin{array}{l}
\displaystyle d\rho^2(\zeta_t^1, \xi_t^1)\le 2\,\langle\dot\gamma(\zeta_t^1, \xi_t^1)(1), b_t(m(\xi_t), \xi^1_t)\rangle \,dt\\
 \\
\displaystyle\hskip3cm -2\,\langle \dot\gamma(\zeta_t^1, \xi_t^1)(0), b_t(\mu_t,\zeta^1_t)\rangle\,dt -\kappa\,\rho(\zeta_t^1, \xi_t^1)^2\,dt.
 \end{array}
\end{equation}
Define similarly to the previous section for a Riemannian manifold $M$ and a map $G : M\times M \to TM$ such that $G(x,y)\in T_{y}M$ : for $z=(x,y)$, $\bar z=(\bar x, \bar y)$ elements of $M\times M$
\begin{equation}
 \label{E4.3.9}
 G[z;\bar z]:=\int_0^1//_{y,\gamma(y,\bar y)(\epsilon)}^{-1}G(\gamma(z,\bar z)(\epsilon))\,d\epsilon~\in~ T_yM.
\end{equation}
Also define 
\begin{equation}
 \label{E4.3.10}
 \Ba_t^{(0)}(z,\bar z):=\left[
 \begin{array}{cc}
  \nabla_yb_t[(y,x);(\bar y,\bar x)]& \nabla_xb_t[(y,x);(\bar y,\bar x)]\\
   \nabla_xb_t[(x,y);(\bar x,\bar y)]& \nabla_yb_t[(x,y);(\bar x,\bar y)]
 \end{array}
 \right],
\end{equation}
\begin{equation}
 \label{E4.3.11}
 \Ba_t^{(1)}(z,\bar z):=\left[
 \begin{array}{cc}
  \nabla(b_t\circ \delta)[ x; \bar x]& 0\\
   0& \nabla(b_t\circ \delta)[ y; \bar y]
 \end{array}
 \right],
\end{equation}
where $\delta : M\to M\times M$, $x\mapsto (x,x)$, and set
\begin{equation}
 \label{E4.3.12}
 \Ba_t(z,\bar z):=\frac1N\,\Ba_t^{(1)}(z,\bar z)+\left(1-\frac1N\right)\Ba_t^{(0)}(z,\bar z)
\end{equation}
Consider the following regularity condition:
  
     {\em $(\Ha_{\Ca}^g)$ : There exists some $\lambda_{\Ca}\in \RR$ such that for any $z,\overline{z}\in M\times M$ and $t\geq 0$ we have 
 \begin{equation}\label{4.3.13}
\Ca_t(z,\overline{z}):= \frac{1}{2}\left[\Ba_{t}(z,\overline{z})+ \Ba_{t}(z,\overline{z})^{\prime}\right]\leq -\lambda_{\Ca}~g_{M\times M}(z)
 \end{equation}}
 \begin{theo}
  \label{T3.3}
  Assume that the Ricci curvatures of $M$ are bounded below by $\kappa\in \RR$ and that the condition $(\Ha_{\Ca}^g)$ is satisfied. Then 
  \begin{equation}
  \label{E4.3.14}
   \EE\left[\rho^2(\zeta_t^1,\xi_t^1)\right]^{1/2}\le \frac{2}{2\lambda_{\Ca}+\kappa}\left(1-e^{-\frac{(2\lambda_{\Ca}+\kappa)t}{2}}\right)~\sqrt{\frac{\beta_t(\mu)}{N}}
  \end{equation}
with the parameter $ \beta_t(\mu)$ defined as in  Theorem~\ref{T4.2}.
 \end{theo}
\begin{remark} The result 
 of Theorem~\ref{T3.3} extends to the case when $\sigma=\sigma_t$ and  $g=g_t$ depend on time, if we replace the bound below of the Ricci curvatures by the assumption that ${\rm Ric}_M-\dot g\ge \kappa g$.
\end{remark}
\proof
 The proof is completely similar to the one of Theorem~\ref{T4.2}, thus it is only sketched. Letting $S_t:= \EE\left[\rho^2(\zeta_t^1,\xi_t^1)\right]$ we arrive at 
\begin{equation}
 \label{E4.3.8bis}
 \partial_tS_t\le (2I_t+2J_t -\kappa\,S_t)\,dt
\end{equation}
where 
\begin{equation}
\label{ItM}
I_t:=\EE\left[\langle\dot\gamma(\zeta_t^1, \xi_t^1)(1), b_t(m(\xi_t), \xi^1_t)\rangle -\langle \dot\gamma(\zeta_t^1, \xi_t^1)(0), b_t(m(\zeta_t), \zeta^1_t)\rangle\,dt\right]
\end{equation}
and 
\begin{equation}
\label{JtM}
J_t:=\EE\left[\langle\dot\gamma(\zeta_t^1, \xi_t^1)(0), b_t(m(\zeta_t), \zeta^1_t)\rangle -\langle \dot\gamma(\zeta_t^1, \xi_t^1)(0), b_t(\mu_t,\zeta^1_t)\rangle\,dt\right]
\end{equation}
which leads to
 \begin{equation}
  \label{E4.3.16}
  \partial_tS_t\le -(2\lambda_{\Ca}+\kappa)S_t+\frac2{\sqrt{N}}\sqrt{S_t}\sqrt{\beta_t(\mu)}
 \end{equation}
 so letting $s_t=\sqrt{S_t}$ we get 
 \begin{equation}
  \label{E4.3.17}
  \partial_ts_t\le -\frac{2\lambda_{\Ca}+\kappa}2s_t+\sqrt{\frac{\beta_t(\mu)}{N}}
 \end{equation}
This ends the proof of (\ref{E4.3.14}).
\cqfd

Let us investigate condition~\eqref{4.3.13} for the Langevin diffusion with drift~\eqref{E3.17}, namely 
\begin{equation}
\label{E4.3.18}
b_t(x,y)=-\nabla U(y) - \nabla (F\circ\rho_x)(y)
\end{equation}
We need the additional assumption $\partial F(0)=0$. In this situation, the computation of $I_t$ in~\eqref{ItM} yields the formula
\begin{equation}
 \label{E4.3.19}
 \begin{split}
 I_t=&-\frac12\EE\left[\nabla^2U^{\oplus 2}\left[(\zeta_t^1,\zeta_t^2);(\xi_t^1,\xi_t^2)\right]\left(\overrightarrow{\zeta_t^1\xi_t^1},\overrightarrow{\zeta_t^2\xi_t^2}\right)\left(\overrightarrow{\zeta_t^1\xi_t^1},\overrightarrow{\zeta_t^2\xi_t^2}\right)\right]\\
 &-\frac12\left(1-\frac1N\right)\EE\left[\nabla^2(F\circ \rho)\left[(\zeta_t^1,\zeta_t^2);(\xi_t^1,\xi_t^2)\right]\left(\overrightarrow{\zeta_t^1\xi_t^1},\overrightarrow{\zeta_t^2\xi_t^2}\right)\left(\overrightarrow{\zeta_t^1\xi_t^1},\overrightarrow{\zeta_t^2\xi_t^2}\right)\right]
 \end{split}
\end{equation}
where we denoted $\overrightarrow{\zeta_t^1\xi_t^1}=\dot\gamma(\zeta_t^1,\xi_t^1)(0)$,
leading to the condition $(\Ha_{\Ca}^g)$: for all $z,\bar z\in M\times M$, 
\begin{equation}
 \label{E4.3.20}
 \nabla^2U^{\oplus 2}[z;\bar z]+\left(1-\frac1N\right)\nabla^2(F\circ \rho)[z;\bar z]\ge \lambda_{\Ca}~g_{M\times M}(z).
\end{equation}
This condition is met for instance when for all $z\in M\times M$, 
\begin{equation}
\label{E4.3.21}
 \nabla^2U^{\oplus 2}(z)+\left(1-\frac1N\right)\nabla^2(F\circ \rho)(z)\ge \lambda_{\Ca}~g_{M\times M}(z).
\end{equation}

\section*{Appendix}
\subsection*{Proof of (\ref{ref-nablax-estimate-0}) and (\ref{ref-nablax-estimate-0-ae})}\label{ref-nablax-estimate-proof}
 After some calculations we check that
\begin{equation}\label{def-C}
 \begin{array}{l}
\displaystyle d \,\left[\nabla X^{\mu}_{s,t}(x) \,\nabla X^{\mu}_{s,t}(x)^{\prime}\right]
=\nabla X^{\mu}_{s,t}(x) ~\Aa_t\left(\phi_{s,t}(\mu),X^{\mu}_{s,t}(x)\right)~\nabla X^{\mu}_{s,t}(x)^{\prime}~dt+dM_{s,t}^{\mu}(x)
\end{array} 
\end{equation} 
with the matrix valued martingale
$$
dM_{s,t}^{\mu}(x):=\sum_{1\leq k\leq r}~\nabla X^{\mu}_{s,t}(x) \left[\nabla_y\sigma_{t,k}\left(\phi_{s,t}(\mu),X^{\mu}_{s,t}(x)\right)+\nabla_y\sigma_{t,k}\left(\phi_{s,t}(\mu),X^{\mu}_{s,t}(x)\right)^{\prime}\right]\nabla X^{\mu}_{s,t}(x)^{\prime}~dW^k_t
$$
and
$$
\Aa_t\left(\mu,y\right):=\nabla_y b_t\left(\mu,y\right)+\nabla_y b_t\left(\mu,y\right)^{\prime}+\sum_{1\leq k \leq r}\nabla_y \sigma_{t,k}\left(\mu,y\right)\nabla_y \sigma_{t,k}\left(\mu,y\right)^{\prime}\leq \int~\mu(dx)~A_t(x,y)
$$
  In the above display, we have used the fact that $ \EE(\Sigma^{\prime})\EE(\Sigma)\leq \EE(\Sigma^{\prime}\Sigma)$, for any random matrix $\Sigma$.
The end of the proof of (\ref{ref-nablax-estimate-0}) and (\ref{ref-nablax-estimate-0-ae}) is now clear.\cqfd

\subsection*{Proof of (\ref{backward-mck-v})}\label{backward-mck-v-proof}

For any time mesh $t_k\leq t_{k+1}$ with $s_0=s$ and $s_n=t$ with $h:=\max{|s_k-s_{k-1}|}$ we have
$$
\PP^{\phi_{s,s_{k-1}}(\mu)}_{s_{k-1},t}=\PP^{\phi_{s,s_{k-1}}(\mu)}_{s_{k-1},s_k}\PP^{\phi_{s_{k-1},s_k}(\phi_{s,s_{k-1}}(\mu))}_{s_{k},t}=
\PP^{\phi_{s,s_{k-1}}(\mu)}_{s_{k-1},s_k}\PP^{\phi_{s,s_k}(\mu)}_{s_{k},t}
$$
Also observe that 
\begin{eqnarray*}
\Delta^{\phi_{s,s_{k-1}}(\mu)}_{s_{k-1},s_k}(x)&:=&X^{\phi_{s,s_{k-1}}(\mu)}_{s_{k-1},s_k}(x)-x\\
&=&\int_{s_{k-1}}^{s_k}~b_u\left(\phi_{s_{k-1},u}(\phi_{s,s_{k-1}}(\mu)),X^{\phi_{s,s_{k-1}}(\mu)}_{s_{k-1},u}(x)\right)~du\\
&&\hskip3cm+\int_{s_{k-1}}^{s_k}~\sigma_u\left(\phi_{s_{k-1},u}(\phi_{s,s_{k-1}}(\mu)),X^{\phi_{s,s_{k-1}}(\mu)}_{s_{k-1},u}(x)\right)~dW_u\\
&=&Y^{\phi_{s,s_{k-1}}(\mu)}_{s_{k-1},s_k}(x)+Z^{\phi_{s,s_{k-1}}(\mu)}_{s_{k-1},s_k}(x)
\end{eqnarray*}
with the random fields
\begin{eqnarray*}
Y^{\mu}_{s_{k-1},s_k}(x)&:=&b_{s_{k}}(\mu,x)~(s_k-s_{k-1})+\sigma_{s_{k}}(\mu,x)~(W_{s_k}-W_{s_{k-1}})\\
Z^{\mu}_{s_{k-1},s_k}(x)&:=&\int_{s_{k-1}}^{s_k}~\left[b_u(\phi_{s_{k-1},u}(\mu),X^{\mu}_{s_{k-1},u}(x))-b_{s_{k}}(\mu,x)\right]~du\\
&&+\int_{s_{k-1}}^{s_k}~\left[
\sigma_u(\phi_{s_{k-1},u}(\mu),X^{\mu}_{s_{k-1},u}(x))-\sigma_{s_{k}}(\mu,x)\right]~dW_u
\end{eqnarray*}

Using elementary manipulations, for any $0\leq h\leq 1$ we check that
$$
       \EE\left(\Vert X^{\mu}_{s,s+h}(x)-x\Vert^{n}\right)^{1/n}\leq~c_n~h~\left[\Vert x\Vert+\mu(\Vert e\Vert^2)^{1/2}\right]\quad \mbox{and}\quad
  \WW_2(\phi_{s,s+h}(\mu),\mu)\leq c~h~\mu(\Vert e\Vert^2)^{1/2}
$$
for some finite constants $c$ and $c_n$.
Recalling that $(t,x,y)\mapsto b_t(x,y)$ and $(t,x,y)\mapsto \sigma_t(x,y)$ are Lipschitz functions we check that
the almost sure convergence 
 $$
 Y^{\mu}_{s_{k-1},s_k}(x)\longrightarrow_{h\rightarrow 0}~0\qquad
  Z^{\mu}_{s_{k-1},s_k}(x)\longrightarrow_{h\rightarrow 0}~0\quad \mbox{and}\quad \Delta^{\mu}_{s_{k-1},s_k}(x)\longrightarrow_{h\rightarrow 0}~0
 $$

Using the Taylor expansion
\begin{eqnarray*}
\PP^{\mu}_{s,t}(f)(x+y)&=&\PP^{\mu}_{s,t}(f)(x)+\tr\left[ \nabla \PP^{\mu}_{s,t}(f)(x)\,y^{\prime}\right]+\frac{1}{2}~\tr\left[\nabla^2 \PP^{\mu}_{s,t}(f)(x) yy^{\prime}\right]\\
&&+\int_0^1(1-\epsilon)~\tr\left[\left(\nabla^2 \PP^{\mu}_{s,t}(f)(x+\epsilon y)-\nabla^2 \PP^{\mu}_{s,t}(f)(x)\right) yy^{\prime}\right]~d\epsilon
\end{eqnarray*}
we check that
$$
\begin{array}{l}
\displaystyle\PP^{\mu}_{s_k,t}(f)\left(x+\Delta^{\mu}_{s_{k-1},s_k}(x)\right)-\PP^{\mu}_{s_k,t}(f)(x)\\
\\
\displaystyle=\tr\left[ \nabla \PP^{\mu}_{s_k,t}(f)(x)\,\Delta^{\mu}_{s_{k-1},s_k}(x)^{\prime}\right]+\frac{1}{2}~\tr\left[\nabla^2 \PP^{\mu}_{s_k,t}(f)(x) 
\Delta^{\mu}_{s_{k-1},s_k}(x)\Delta^{\mu}_{s_{k-1},s_k}(x)^{\prime}
\right]\\
\\
\displaystyle+\int_0^1(1-\epsilon)~\tr\left[\left(\nabla^2 \PP^{\mu}_{s_k,t}(f)(x+\epsilon \Delta^{\mu}_{s_{k-1},s_k}(x))-\nabla^2 \PP^{\mu}_{s_k,t}(f)(x)\right) \Delta^{\mu}_{s_{k-1},s_k}(x)\Delta^{\mu}_{s_{k-1},s_k}(x)^{\prime}\right]~d\epsilon
\end{array}
$$
Rearranging the terms we find  that
$$
\begin{array}{l}
\displaystyle\PP^{\mu}_{s_k,t}(f)\left(x+\Delta^{\mu}_{s_{k-1},s_k}(x)\right)-\PP^{\mu}_{s_k,t}(f)(x)\\
\\
\displaystyle=\tr\left[ \nabla \PP^{\mu}_{s_k,t}(f)(x)\,Y^{\mu}_{s_{k-1},s_k}(x)^{\prime}\right]+\frac{1}{2}~\tr\left[\nabla^2 \PP^{\mu}_{s_k,t}(f)(x) Y^{\mu}_{s_{k-1},s_k}(x)Y^{\mu}_{s_{k-1},s_k}(x)^{\prime}\right]+R^{\mu}_{s_k,t}(f)(x)
\end{array}
$$
with the remainder term
$$
\begin{array}{l}
\displaystyle R^{\mu}_{s_k,t}(f)(x):=\tr\left[ \nabla \PP^{\mu}_{s_k,t}(f)(x)\,Z^{\mu}_{s_{k-1},s_k}(x)^{\prime}\right]+\frac{1}{2}~\tr\left[\nabla^2 \PP^{\mu}_{s_k,t}(f)(x) Z^{\mu}_{s_{k-1},s_k}(x)Z^{\mu}_{s_{k-1},s_k}(x)^{\prime}\right]\\
\\
\displaystyle+\int_0^1(1-\epsilon)~\tr\left[\left(\nabla^2 \PP^{\mu}_{s_k,t}(f)(x+\epsilon \Delta^{\mu}_{s_{k-1},s_k}(x))-\nabla^2 \PP^{\mu}_{s_k,t}(f)(x)\right) \Delta^{\mu}_{s_{k-1},s_k}(x)\Delta^{\mu}_{s_{k-1},s_k}(x)^{\prime}\right]~d\epsilon
\end{array}
$$
This yields the decomposition
$$
\begin{array}{l}
\displaystyle\PP^{\mu}_{s_k,t}(f)\left(x+\Delta^{\mu}_{s_{k-1},s_k}(x)\right)-\PP^{\mu}_{s_k,t}(f)(x)\\
\\
\displaystyle=\left\{\tr\left[ \nabla \PP^{\mu}_{s_k,t}(f)(x)\,b_{s_k}(\mu,x)^{\prime}~\right]+\frac{1}{2}~\tr\left[\nabla^2 \PP^{\mu}_{s_k,t}(f)(x) ~\sigma_{s_k}(\mu,x)~\sigma_{s_k}(\mu,x)^{\prime}\right]\right\}~(s_k-s_{k-1})\\
\\
\hskip3cm+\tr\left[ \nabla \PP^{\mu}_{s_k,t}(f)(x)\,~(W_{s_k}-W_{s_{k-1}})^{\prime}~\sigma_{s_k}(\mu,x)^{\prime}\right]+\overline{R}^{\mu}_{s_k,t}(f)(x)
\end{array}
$$
with the remainder term
$$
\begin{array}{l}
\displaystyle \overline{R}^{\mu}_{s_k,t}(f)(x):=R^{\mu}_{s_k,t}(f)(x)+\frac{1}{2}~\tr\left[\nabla^2 \PP^{\mu}_{s_k,t}(f)(x) ~b_{s_k}(\mu,x)\,(W_{s_k}-W_{s_{k-1}})^{\prime}\sigma_{s_k}(\mu,x)^{\prime}\right](s_k-s_{k-1})\\
\\
\hskip1cm\displaystyle+\frac{1}{2}~\tr\left[\nabla^2 \PP^{\mu}_{s_k,t}(f)(x) ~\sigma_{s_k}(\mu,x)~(W_{s_k}-W_{s_{k-1}})~b_{s_k}(\mu,x)^{\prime}\right]~(s_k-s_{k-1})\\\
\\
\hskip1.5cm\displaystyle+\frac{1}{2}~\tr\left[\nabla^2 \PP^{\mu}_{s_k,t}(f)(x) ~b_{s_k}(\mu,x))b_{s_k}(\mu,x)^{\prime}~\right]~(s_k-s_{k-1})^2\\
\\
\hskip1.5cm+\displaystyle\frac{1}{2}~\tr\left[\nabla^2 \PP^{\mu}_{s_k,t}(f)(x)~ \sigma_{s_k}(\mu,x)~\left[(W_{s_k}-W_{s_{k-1}})(W_{s_k}-W_{s_{k-1}})^{\prime}-(s_k-s_{k-1})~I\right]\sigma_{s_k}(\mu,x)^{\prime}\right]
\end{array}
$$
On the other hand, we have
\begin{eqnarray*}
\PP^{\mu}_{s,t}(f)(x)-f(x)&=&\sum_{1\leq k\leq n}\left[\PP^{\phi_{s,s_{k-1}}(\mu)}_{s_{k-1},t}(f)(x)-\PP^{\phi_{s,s_{k}}(\mu)}_{s_{k},t}(f)(x)\right]\\
&=&\sum_{1\leq k\leq n}\left[\PP^{\phi_{s,s_{k}}(\mu)}_{s_{k},t}(f)\left(
x+\Delta^{\phi_{s,s_{k-1}}(\mu)}_{s_{k-1},s_k}(x)\right)-\PP^{\phi_{s,s_{k}}(\mu)}_{s_k,t}(f)(x)\right]
\end{eqnarray*}
This implies that
$$
\begin{array}{l}
\displaystyle\PP^{\mu}_{s,t}(f)(x)-f(x)=\sum_{1\leq k\leq n}\left\{\tr\left[ \nabla \PP^{\phi_{s,s_{k}}(\mu)}_{s_k,t}(f)(x)\,b_{s_k}(\phi_{s,s_{k}}(\mu),x)^{\prime}~\right]\right.\\
\\
\hskip3cm\displaystyle\left.+\frac{1}{2}~\tr\left[\nabla^2 \PP^{\phi_{s,s_{k}}(\mu)}_{s_k,t}(f)(x) ~\sigma_{s_k}(\phi_{s,s_{k}}(\mu),x)~\sigma_{s_k}(\phi_{s,s_{k}}(\mu),x)^{\prime}\right]\right\}~(s_k-s_{k-1})\\
\\
\hskip3cm\displaystyle+\tr\left[ \nabla \PP^{\phi_{s,s_{k}}(\mu)}_{s_k,t}(f)(x)\,~(W_{s_k}-W_{s_{k-1}})^{\prime}~\sigma_{s_k}(\phi_{s,s_{k}}(\mu),x)^{\prime}\right]+\sum_{1\leq k\leq n}\overline{R}^{\phi_{s,s_{k}}(\mu)}_{s_k,t}(f)(x)
\end{array}
$$
We end the proof of (\ref{backward-mck-v}) by letting the time step $h\rightarrow 0$.\cqfd
\subsection*{Proof of theorem~\ref{theo-Aa}}\label{c-rule-proof}
Observe that
\begin{eqnarray*}
\partial_{z^l_i}\left[\Ha_1^{j,k}\circ\Ha_2\right](z)&=&\sum_{1\leq m\leq d}\sum_{1\leq n\leq N}\left(\partial_{z^m_n}\Ha_1^{j,k}\right)(\Ha_2(z))~\partial_{z^l_i}\Ha_2^{n,m}(z)\\
&=&\sum_{1\leq m\leq d}\sum_{1\leq n\leq N}~\left[\nabla_{z_i}\Ha_2^n\right]_{l,m}(z)~\left[\nabla_{z_n}\Ha_1^j\right]_{m,k}(\Ha_2(z))
\end{eqnarray*}
This implies that
$$
\displaystyle\left[\nabla\left(\Ha_1\circ\Ha_2\right)(z)\right]_{i,j}=\nabla_{z_i}\left[\Ha^{j}_1(\Ha_2(z))\right]=\sum_{1\leq n\leq N}\left[\nabla\Ha_2(z)\right]_{i,n}~\left[\left(\nabla \Ha_1\right)(\Ha_2(z))\right]_{n,j}
$$
This ends the proof of  (\ref{c-rule}). The proof of (\ref{Ha-ref}) and (\ref{Ha-ref-2}) 
come from the formula
$$
 d\,\nabla\xi_{s,t}(z)\nabla\xi_{s,t}(z)^{\prime}=\nabla\xi_{s,t}(z)\Aa_t\left(\xi_{s,t}(z)\right)
\nabla\xi_{s,t}(z)^{\prime}dt+d\Ma_t(z)
$$
with the martingale
$$\Ma_t(z)_{i,k}=\sum_{1\leq j\leq N}~\sum_{1\leq \alpha\leq d}\Ma^{j,\alpha}_t(z)_{i,k}$$
defined in terms of the diffusion processes
$$
\begin{array}{l}
\displaystyle d\Ma^{j,\alpha}_t(z)_{i,k}
:=
 \left[\nabla\xi_{s,t}(z)\nabla\Ga_{t,\alpha}\left(\xi_{s,t}(z)\right)\right]_{i,j}~ dW^{j,\alpha}_t ~
\nabla\xi_{s,t}(z)^{\prime}_{j,k}\\
 \\
 \hskip3cm+
~
\nabla\xi_{s,t}(z)_{i,j}~ dW^{j,\alpha}_t ~
\left(\left[\nabla\xi_{s,t}(z)\right]\left[\nabla\Ga_{t,\alpha}\right]\left(\xi_{s,t}(z)\right)\right)^{\prime}_{j,k}
\end{array}
$$
The end of the proof of (\ref{Ha-ref}) and (\ref{Ha-ref-2}) follows the same lines of arguments as the proof of (\ref{ref-nablax-estimate-0}) and (\ref{ref-nablax-estimate-0-ae}), thus it is skipped. This ends the proof of the theorem.
\cqfd


\begin{thebibliography}{99}
 \bibitem{ACT10}
 M. Arnaudon, A. Coulibaly, A. Thalmaier,
Horizontal diffusions in $C^1$ path space.
S\'eminaire de Probabilit\'es XLIII, 73--94, Lecture Notes in Mathematics 2006, Springer (2010)
 
 
 \bibitem{ATW06}
 M. Arnaudon, A. Thalmaier and F.Y. Wang,
 Harnack inequality and heat kernel estimates on manifolds with curvature unbounded below,
 Bull. Sci. Math. 130 (2006) 223--233
 
 
 \bibitem{bakry}
D. Bakry, M. Emery. Diffusions hypercontractives. In S\'eminaire de probabilit\'es, XIX, 1983/84, pp. 177--206. 
Lect. Notes in Math. 1123, Springer (1985).
 \bibitem{benachour-1}
S. Benachour, B. Roynette, D. Talay, P. Vallois. Nonlinear self-stabilizing processes, part I. Existence, invariant probability, propagation of chaos. Stochastic processes and their applications, vol. 75, no. 2, pp. 173--201 (1998).


\bibitem{benachour-2}
S. Benachour, B. Roynette, P. Vallois. Nonlinear self-stabilizing processes, part II: Convergence to invariant probability. Stochastic processes and their applications, vol. 75, no.2, pp. 203--224 (1998).

 
 \bibitem{bene-1}
 D. Benedetto, E. Caglioti, M. Pulvirenti. A kinetic equation for granular media. RAIRO Mod\`el. Math. Anal. Num\'er. vol. 31, no. 5, pp. 615--641 (1997).
 
  \bibitem{bene-2}
   D. Benedetto, E. Caglioti, E., Carrillo, M. Pulvirenti. A non-Maxwellian steady distribution for one-dimensional granular media. J. Statist. Phys.vol.  91, pp. 979--990 (1998).
   
   \bibitem{gbg-13}
F. Bolley, I. Gentil, A. Guillin. Uniform convergence to equilibrium for granular media.
Archive for Rational Mechanics and Analysis, vol. 208, no. 2, pp. 429--445 (2013).
   
   \bibitem{calvez}
V. Calvez, Vincent, J. Carrillo. Refined asymptotics for the subcritical Keller-Segel system and related functional inequalities.
Proceedings of the American Mathematical Society, vol. 140, no. 10, pp. 3515--3530 (2012).

\bibitem{cmcv-2006}
J. A. Carrillo, R. J. McCann, and C. Villani. Contractions in the 2-Wasserstein length space and thermalization of granular media. Arch. Rational Mech. Anal., vol. 179, pp. 217--263, (2006).

\bibitem{c-toscani}
J. A. Carrillo and G. Toscani. Wasserstein metric and large-time assymptotics of nonlinear diffusion equations. In New trends in math. physics. World Sci., Singapore (2005).

\bibitem{cattiaux}
P. Cattiaux, A. Guillin, and F. Malrieu. Probabilistic approach for granular media equations in the non uniformly convex case. Prob. Theor. Rel. Fields, vol. 140, no. 1-2, pp. 19--40 (2008).

\bibitem{carverhill}
A.P. Carverhill and K.D. Elworthy, Flows of Stochastic Dynamical Systems: The Functional Analytic Approach. Z. Wahrs vol. 65, pp. 245-267 (1983).
 
   \bibitem{coppel}
 W. A. Coppel. Stability and asymptotic behavior of differential equations. D. C. Heath and Co., Boston, Mass. (1965).
 
 \bibitem{cordero}
D. Cordero-Erausquin, W. Gangbo, and C. Houdr\'e. Inequalities for generalized entropy and optimal transportation. In Recent Advances in the Theory and Applications of Mass Transport, Contemp. Math. 353. A. M. S., Providence (2004).
 
  \bibitem{daprato-2}
 G. Da Prato, 
J.L.  Menaldi, L.  Tubaro. Some results of backward Ito formula. Stochastic analysis and applications, vol. 25, no. 3, pp. 679--703 (2007).

 
 \bibitem{tugaut-hong}
 M.H. Duong, J. Tugaut. The Vlasov-Fokker-Planck equation in non-convex landscapes: convergence to equilibrium. Electronic Communications in Probability, vol. 23 (2018).
 
 
\bibitem{graham}
C. Graham. McKean-Vlasov, Ito-Skorohod equations and nonlinear diffusions with discrete jumps. Stochastic Processes and their Applications.
vol. 40, pp. 69-82 (1992).


\bibitem{h-t-2010}
S. Herrmann, J. Tugaut. Non-uniqueness of stationary measures for self-stabilizing processes.
 Stochastic Processes and their Applications, vol.120, no. 7, pp. 1215--1246 (2010).
 
 
 \bibitem{2-h-t-2010}
S. Herrmann, J. Tugaut. Stationary measures for self-stabilizing processes: asymptotic analysis in the small noise limit.  Electron. J. Probab., vol. 15, pp. 2087--2116  (2010).


 \bibitem{3-h-t-2010}
S. Herrmann, J. Tugaut. Self-stabilizing processes: uniqueness problem for stationary measures and convergence rate in the small-noise limit. ESAIM: Probability and Statistics, vol. 16, pp. 277--305 (2012).

\bibitem{huang}
X. Huang, M. R\"ockner, F.Y. Wang. . Nonlinear Fokker-Planck equations for probability measures on path space and path-distribution dependent SDEs. arXiv preprint arXiv:1709.00556 (2017).



\bibitem{khalil}
H.K. Khalil, Nonlinear Systems. Prentice Hall (2002).

\bibitem{lyapunov}
 A.M. Lyapunov, The general problem of the stability of motion, CRC Press (1992).
 

\bibitem{malrieu}
F. Malrieu. Logarithmic Sobolev inequalities for some nonlinear PDE's. Stochastic Process. Appl. vol. 95,  no. 1, 109?132 (2001).

\bibitem{malrieu-2}
F.  Malrieu. Convergence to equilibrium for granular media equations and their Euler schemes. Ann. Appl. Probab., vol. 13, no. 2, pp. 540--560 (2003).

\bibitem{martin}
R.H Martin Jr. Bounds for solutions of a class of nonlinear differential equations. Journal of Differential Equations, vol. 8, no. 3, pp.416 --430 (1970).


 

 \bibitem{mckean-1}
 H. P. McKean. A class of Markov processes associated with nonlinear parabolic equations. Proc. Nat. Acad. Sci. U.S.A., vol. 56, pp. 1907--1911 (1966).

 \bibitem{mckean-2}
 H. P. McKean. Propagation of chaos for a class of non-linear parabolic equations. In Stochastic Differential Equations, Lecture Series in Differential Equations, Session 7, Catholic Univ. (1967), pp. 41--57. Air Force Office Sci. Res., Arlington, Va., (1967).
 

\bibitem{norris}
J.R. Norris. Simplified Malliavin calculus. S\'eminaire de Probabilist\'es, vol. 20, pp. 101-130 (1986).

\bibitem{otto}
F. Otto. The geometry of dissipative evolution equations: the porous medium equation. Comm. Partial Differential Equations, vol. 26, pp. 101--174 (2001).

\bibitem{otto-2}
F. Otto, V. Villani. Generalization of an inequality by Talagrand, and links with the logarithmic Sobolev inequality. J. Funct. Anal. vol. 
173, pp. 361--400 (2000).

\bibitem{tamura}
Y. Tamura. On asymptotic behaviors of the solution of a nonlinear diffusion equation. J. Fac. Sci. Univ. Tokyo Sect. IA Math., vol.31, no. 1, pp. 195--221 (1984).

\bibitem{tamura-2}
Y. Tamura. Free energy and the convergence of distributions of diffusion processes of McKean type. J. Fac. Sci. Univ. Tokyo Sect. IA Math., vol. 34, no. 2, pp. 443--484, (1987).
\bibitem{toscani}
G. Toscani. One-dimensional kinetic models of granular flows. RAIRO Mod\`el. Math. Anal. Num\'er. no. 34, no. 6, pp. 1277--1291 (2000).
\bibitem{tugaut-1}
J. Tugaut. Convergence to the equilibria for self-stabilizing processes in double well landscape.  Annals of Probability, vol. 41, no. 3A, pp. 1427--1460 (2013).

\bibitem{tugaut-2}
J. Tugaut. Self-stabilizing processes in multi-wells landscape in $\RR^d$ - Invariant probabilities.  J Theor Probab. vol. 27, no. 1, pp. 27--57 (2014).


\bibitem{tugaut-3}
J. Tugaut. Convergence in Wasserstein distance for self-stabilizing diffusion evolving in a double-well landscape. Comptes Rendus Mathematiques,  vol. 356, no. 6,  pp. 657--660 (2018).

\bibitem{tugaut-4}
J. Tugaut.  Exit problem of McKean-Vlasov diffusions in convex landscapes. Electronic Journal of Probability, vol. 17, no. 76 , pp. 1-26 (2012).

\bibitem{tugaut-5}
J. Tugaut. 
A simple proof of a Kramers' type law for self-stabilizing diffusions. Electronic Communications in Probability, vol. 21 (2016).


\bibitem{tugaut-6}
J. Tugaut. 
McKean-Vlasov diffusions: from the asynchronization to the synchronization. CR Math. Acad. Sci, vol. 349, no. 17-18, pp. 983-986  (2011).

\bibitem{villani}
C. Villani. A survey of mathematical topics in the collisional kinetic theory of gases. Handbook of mathematical fluid dynamics, vol. 1, no 71--305, p. 3--8 (2002).

\bibitem{WangFY14}
Feng-Yu Wang,
Analysis for diffusion processes on Riemannian manifolds, Advanced Series on Statistical Science and Applied Probability, Vol. 18, World Scientific (2014).
\end{thebibliography}
\end{document}